\documentclass[11pt,twoside,a4paper]{article}
\usepackage{graphicx}
\usepackage{amsfonts}
\usepackage{bbm}
\usepackage[leqno]{amsmath}%%%%%%%%%%%%%%%%%%%%%%%%%%%%%%%%%%%% number the equations at the lefthand side%%%%%%%%%%%%%%%%%%%%%%%%%%%%%%%%%%%55
\usepackage{mathrsfs}
\usepackage{fancyhdr}
 \usepackage{mathrsfs,amsmath,amssymb,amsthm}
\usepackage[dvips]{color}
 \usepackage{amssymb}
 \usepackage{amsbsy}
 \usepackage[dvips]{color}
\usepackage{titlesec}
\titlespacing*{\section}{0pt}{0.5\baselineskip}{0.5\baselineskip}
\titleformat*{\subsubsection}{\it}
\titlespacing*{\subsection}{0pt}{0.3\baselineskip}{0.3\baselineskip}
\titlespacing*{\subsubsection}{0pt}{0.3\baselineskip}{0.3\baselineskip}

\allowdisplaybreaks

 \theoremstyle{definition}
\theoremstyle{remark}  \newtheorem{remark}{\noindent\mbox{Remark}}
 \theoremstyle{plain}
 \theoremstyle{plain}\newtheorem{lemma}{\noindent\mbox{Lemma}}
\theoremstyle{plain} \newtheorem{theorem}{\noindent\mbox{Theorem}}
 \theoremstyle{plain}\newtheorem{proposition}{\noindent\mbox{Proposition}}
 \theoremstyle{plain}
\theoremstyle{definition} 

 \def\proof{\noindent{\it Proof.~~}}
 \def\qed{\hfill$\Box$\medskip}
 \def\rto{\rightarrow\infty}
\def\z{\left}
\def\y{\right}
 \def\no{\nonumber}
 \def\mb{\mathbf}

\begin{document}

 \title{{On extinction time distribution of a 2-type linear-fractional branching process in a varying environment with asymptotically constant mean matrices}}                %%%   the Fund which you are supported by  %%%

\author{Hua-Ming \uppercase{Wang}\footnote{Email:hmking@ahnu.edu.cn; School of Mathematics and Statistics, Anhui Normal University, Wuhu, 241003, China }}

\date{}
\maketitle%

\vspace{-1cm}

\begin{center}
\begin{minipage}[c]{12cm}
\begin{center}\textbf{Abstract}\quad \end{center}
In this paper we study a 2-type linear-fractional branching process in varying environment with asymptotically constant mean matrices.  Let $\nu$ be the extinction time and for $k\ge1$ let $M_k$ be the mean matrix of offspring distribution of individuals of the $(k-1)$-th generation. Under certain conditions, we show that $P(\nu=n)$ and $P(\nu>n)$ are asymptotically equivalent to some functions of products of spectral radii of the mean matrices. This paper complements
 a former result [{\it arXiv:} 2007.07840] which requires in addition a condition $\forall k\ge1,\rm{det}(M_k)<-\varepsilon$ for some $\varepsilon>0.$  Such a condition excludes a large class of mean matrices. As byproducts, we also get some results on asymptotics of products of nonhomogeneous matrices which have their own interests.

\vspace{0.2cm}

\vspace{0.2cm}

\end{minipage}
\end{center}

\section{Introduction}

{\bf 1.1 Background and motivation.} It is known that many new phenomena, which homogeneous Galton-Watson processes do not possess, arise when considering the branching processes in varying environments (BPVE hereafter). Thus BPVE has been extensively studied by many authors. For details of the single-type case, we refer the reader to \cite{bp}, \cite{fuj}, \cite{j},  \cite{ker}, \cite{va}, \cite{lin}, \cite{ms}  and references therein.

Compared with the single-type case, the situation of the multitype setting is less satisfying. The   convergence of the normalized population size and criteria for almost sure extinction can be found in
 \cite{bcn}, \cite{cw}, \cite{jon} and a recent article \cite{dhkp}.
   We notice also that in \cite{dhkp}, Dolgopyat et al. gave also the asymptotics of the survival probabilities. Their proof relies on a generalization of the Perron-Frobenius theorem suitable for studying the product of nonhomogeneous nonnegative matrices. Let $\nu$ be the extinction time and $M_k$ be the mean matrix of offspring distribution of individuals of the $(k-1)$-th generation. They showed that \begin{align}\label{sv}
  \frac{1}{C}\Big(\sum_{k=1}^n\lambda_1^{-1}\cdots\lambda_{n}^{-1}\Big)^{-1}<P(\nu>n)<C\Big(\sum_{k=1}^n\lambda_1^{-1}\cdots\lambda_{n}^{-1}\Big)^{-1}, \end{align} where $0<C<\infty$ is a certain constant and for $k\ge1,$ $\lambda_k$ is a number associated to the generalized Perron-Frobenius theorem and depends on the mean matrices $M_n,n\ge k$ of the branching process, see \cite[Proposition 2.1 and Lemma 2.2]{dhkp}.
 Although \eqref{sv} provides a lower and an upper bound for the survival probability, it does not give the asymptotical equivalence of the survival probability and furthermore, those number $\lambda_k,k\ge1$ are hard to compute explicitly. For this consideration, Wang and Yao \cite{wy} considered a 2-type linear-fractional branching process with asymptotically constant mean matrices. Instead of those number $\lambda_k,$ they used the spectral radii of the mean matrices of the offspring distributions, which can be explicitly computed.
For this special setting, they give not only the asymptotical equivalence of $P(\nu>n)$ but also that of $P(\nu=n).$  But, they need a condition
\begin{align}\label{p}
  \forall k\ge1, M_k(12)M_{k}(21)>M_k(11)M_k(22)+\varepsilon
\end{align} for some $\varepsilon>0,$ which excludes a large class of mean matrices.  They added such a restriction because their proof depends on some result on the asymptotical equivalence between the elements of products of {\it positive} matrices and the products of spectral radii of those matrices, which was studied in \cite{hs20}, and some delicate analysis of the tails and the critical tails of {\it positive} continued fractions.
Without \eqref{p}, one has to deal with product of matrices with some {\it negative elements} and also continued fractions with {\it negative coefficients}.

In this paper, by generalizing the  result of product of matrices in \cite{hs20} and analyzing the tails and approximants of continued fractions with negative coefficients, we remove the assumption \eqref{p}
to characterize asymptotics of the tail probability and the probability that the extinction time equals $n$ as $n\rto.$
%
%Currently, for multitype BPVE with general offspring distribution, we have no idea how to solve the above questions. In this paper, we consider  two-type linear-fractional BPVEs with asymptotically constant mean matrices.
%For the linear-fractional setting, the distribution of the extinction time can be written in terms of sum of product $\prod_{k=1}^n M_k$ of the mean matrices.
%But $\prod_{k=1}^n M_k$ is hard to be estimated. So we construct some new matrix $A_k$ which may depend on $M_k$ and $M_{k+1}.$ The product $\prod_{k=1}^nA_k$ is related to the approximants of a continued fraction. Thus by some delicate analysis of the continued fractions and the product $\prod_{k=1}^n A_k,$  we can express the asymptotics of extinction time distribution in terms of the product of the spectral radii of $M_k,$ which can be computed explicitly.
%
%We also construct examples for which the mass of the extinction time distribution at $n$ decays with various speeds, for example, $\frac{c}{n(\log n)^2},$ $\frac{c}{n^\beta},\beta>1$ et al. This observation complements  Fujimagari \cite{fuj}, which studied the single-type counterpart.

\noindent{\bf 1.2 Model and main results.} Suppose that $M_k,k\ge1$ is a sequence of  nonnegative 2-by-2 matrices and $\gamma_k=(\gamma_k^{(1)},\gamma_k^{(2)}), k\ge 1$ is a sequence of nonnegative  row vectors. To avoid the degenerate case, we require that $\forall k\ge1,$ all elements of $M_kM_{k+1}$ are strictly positive and $\gamma_k\ne \mb 0.$
For $\mathbf s=(s_1,s_2)^t\in [0,1]^2$ and $k\ge1,$ let
\begin{align*}
 \mathbf{f}_{k}(\mathbf{s})=(f_{k}^{(1)}(\mathbf{s}),f_{k}^{(2)}(\mathbf{s}))^{t}=\mb1-\frac{M_{k}(\mb1-\mathbf{s})}{1+\gamma_{k}(\mb1-\mathbf{s})}
\end{align*}
which is known as the probability generating function of a linear-fractional distribution. Here and in what follows, $\mb v^t$ denotes the transpose of a vector $\mb v$ and $\mb 1=(\mb e_1+\mb e_2)^t=(1,1)^t,$ with $\mb e_1=(1,0),\mb e_2=(0,1).$

Suppose that $Z_n=(Z_{n,1},Z_{n,2}),n\ge0$ is a stochastic process such that
\begin{align*}
  E\z(\mathbf s^{Z_n}\big|Z_0,...,Z_{n-1}\y)=\z[\mathbf f_{n}(\mb s)\y]^{Z_{n-1}}, n\ge1,
\end{align*}
where $\z[\mathbf f_{n}(s)\y]^{Z_{n-1}}:=\z[f_n^{(1)}(\mb s)\y]^{Z_{n,1}}\z[f_n^{(2)}(\mb s)\y]^{Z_{n,2}}.$
We call the process $Z_n,n\ge0$ a two-type linear-fractional branching process in a varying environment. Matrices $M_k,k\ge1$ are usually referred to as the mean matrices of the branching process.
Denote by $$\nu=\min\{n: Z_{n}=\mathbf 0\}$$ the extinction time of $\{Z_n\}$ which we concern.

 Throughout, we assume $b_k,d_k>0,$  $a_k,\theta_k\ge0,$ $a_k+\theta_k>0,$ $k\ge1$  and put
\begin{align}\label{mg} M_k:=\left( \begin{array}{cc}
  a_k & b_k \\
  d_k &\theta_k \\
 \end{array}\right),\gamma_k:=\mb e_1M_k,\forall k\ge1.\end{align}

 We introduce the following conditions on the number $a_k,b_k,d_k,\theta_k,k\ge1.$

\noindent{\bf(B1)} Suppose that $b,d>0,$  $ a,\theta\ge0$ are some numbers  such that $a+\theta>0,$
$ a_k\rightarrow a,   b_k\rightarrow b ,  d_k\rightarrow d,\theta_k\rightarrow\theta$ as $k\rto$
 and assume further that
\begin{align}
  \sum_{k=2}^\infty|a_k-a_{k-1}|+| b_k- b_{k-1}|+| d_k- d_{k-1}|+|\theta_{k}-\theta_{k-1}|<\infty.\no
\end{align}

Suppose now condition (B1) holds and for $k\ge1,$ set \begin{align}\label{dta} A_k:=\left( \begin{array}{cc}
  \tilde a_k & \tilde b_k \\
  \tilde d_k &0 \\
 \end{array}\right) \text{ with }
\tilde a_k=a_k+\frac{b_k\theta_{k+1}}{ b_{k+1}}, \tilde b_k= b_k,\tilde  d_k=d_k-\frac{a_k\theta_k}{b_k}.\end{align}
Letting $\Lambda_k=\left(
                 \begin{array}{cc}
                   1 & 0 \\
                   \theta_k/b_k & 1 \\
                 \end{array}
               \right),k\ge1,$ then for $n\ge k\ge1,$ we have
\begin{align}\label{am}
  A_k=\Lambda_k^{-1}M_k\Lambda_{k+1}\text{ and }\mb e_1\prod_{i=k}^n M_i\mb 1=\mb e_1 \prod_{i=k}^n A_i(1,1-\theta_{n+1}/b_{n+1})^t.
\end{align}

For linear-fractional setting, the distribution of $\nu$ can be formulated explicitly by $M_k,$ see \eqref{pngn} and \eqref{pnen} below. However, the elements of $\prod_{i=k}^n M_i$ are hard to compute and evaluate whereas those of $\prod_{i=k}^n A_i$ are workable because they have some correspondence with continued fractions due to the special structure of the matrices $A_i,i\ge1.$ Therefore, instead of $M_k,$ we will work with $A_k$ below.

We need in addition the following conditions which are mutually exclusive.

\noindent{\bf(B2)$_{a}$} $\exists k_0>0,$ such that $\frac{\tilde a_k}{\tilde b_k}=\frac{\tilde a_{k+1}}{\tilde b_{k+1}},
\ \frac{\tilde d_k}{\tilde b_k}\ne\frac{\tilde d_{k+1}}{\tilde b_{k+1}},\ \forall k\ge k_0$ and
$$\lim_{k\rto}\frac{\tilde d_{k+2}/\tilde b_{k+2}-\tilde d_{k+1}/\tilde b_{k+1}}{\tilde d_{k+1}/\tilde b_{k+1}-\tilde d_{k}/\tilde b_{k}}\text{ exists.}$$

 \noindent{\bf(B2)$_{b}$} $\exists k_0>0,$ such that $\frac{\tilde a_k}{\tilde b_k}\ne\frac{\tilde a_{k+1}}{\tilde b_{k+1}},
\ \frac{\tilde d_k}{\tilde b_k}=\frac{\tilde d_{k+1}}{\tilde b_{k+1}},\ \forall k\ge k_0$ and
$$\lim_{k\rto}\frac{\tilde a_{k+2}/\tilde b_{k+2}-\tilde a_{k+1}/\tilde b_{k+1}}{\tilde a_{k+1}/\tilde b_{k+1}-\tilde a_{k}/\tilde b_{k}}\text{ exists.}$$

\noindent{\bf(B2)$_{c}$} $\exists k_0>0,$ such that $\frac{\tilde a_k}{\tilde b_k}\ne\frac{\tilde a_{k+1}}{\tilde b_{k+1}},
\ \frac{\tilde d_k}{\tilde b_k}\ne\frac{\tilde d_{k+1}}{\tilde b_{k+1}},\ \forall k\ge k_0$ and
   $$\tau:=\lim_{k\rto}\frac{\tilde d_{k+1}/\tilde b_{k+1}-\tilde d_{k}/\tilde b_k}{\tilde a_{k+1}/\tilde b_{k+1}-\tilde a_{k}/\tilde b_k}\ne \frac{-(a+\theta)\pm \sqrt{(a+\theta)^2+4(bd-a\theta)}}{2b}  $$ exists as a finite or infinite number. In addition, if $\tau$ is finite,  assume further that $\lim_{k\rto}\frac{\tilde a_{k+2}/\tilde b_{k+2}-\tilde a_{k+1}/\tilde b_{k+1}}{\tilde a_{k+1}/\tilde b_{k+1}-\tilde a_{k}/\tilde b_{k}}
$ exists. Otherwise, if $\tau=\infty,$ assume  further that $\lim_{k\rto}\frac{\tilde d_{k+2}/\tilde b_{k+2}-\tilde d_{k+1}/\tilde b_{k+1}}{\tilde d_{k+1}/\tilde b_{k+1}-\tilde d_{k}/\tilde b_{k}}$ exists.

The conditions (B2)$_a,$ (B2)$_b$ and  (B2)$_c$ look a bit complicated. But they are not so difficult to fulfilled. We refer the reader to \cite[Lemma 1]{wy} for examples for which (B1) and one of (B2)$_a,$ (B2)$_b$ and  (B2)$_c$ hold.
%
%\begin{lemma}\label{egc}
%  Suppose that  $b,d>0$ and $a,\theta\ge0$ are numbers which are not all equal, satisfying $a+\theta>0,$ $\tau:=\frac{b(b+d-a-\theta)+2(a\theta-bd)}{b(2b-a-\theta)}\ne \frac{-(a+\theta)\pm\sqrt{(a+\theta)^2+4(bd-a\theta)}}{2b}$ and $(b-a)(b-\theta)\ge0.$ Let $r_n,n\ge 1$ be strictly positive numbers such that $\lim_{n\rto}r_n=0$ and $\lim_{n\rto}\frac{r_n-r_{n+1}}{r_n^{2}}=c$ for some number $0<c<\infty.$ Set $a_k=a+r_k,$ $b_k=b+r_k,$ $d_k=d+r_k$ and $\theta_k=\theta+r_k,$ $k\ge1.$ Then $a_k,b_k,d_k,\theta_k,k\ge1$ satisfies condition {(B1)} and one of the conditions {(B2)$_{ a},$ (B2)$_{ b}$} and {(B2)$_{ c}.$}
%\end{lemma}

Note that under condition (B1), we have \begin{equation}\no
  \lim_{k\rto}M_k=M:=\z(\begin{array}{cc}
                          a & b \\
                          d & \theta
                        \end{array}
  \y),\ \lim_{k\rto} A_k=A:=\z(\begin{array}{cc}
                          a+\theta & b \\
                          d-a\theta/b & 0
                        \end{array}
  \y)
\end{equation}
whose eigenvalues are  \begin{align}
  \varrho:=\varrho(M)&=\varrho(A)=\frac{a+\theta+\sqrt{(a+\theta)^2+4(bd- a\theta)}}{2},\label{r}\\
\varrho_1:=\varrho_1(M)&=\varrho_1(A)=\frac{a+\theta-\sqrt{(a+\theta)^2+4(bd-a\theta)}}{2}.\label{r1}
\end{align}
It is clear that $|\varrho_1|<\varrho.$ In the literature,  $\varrho(B):=\sup\{|\lambda|: \lambda \text{ is an eigenvalue of }B \}$ is usually referred to as the spectral radius of a matrix $B.$
In what follows, we always denote by $\varrho(B)$  the spectral radius of a matrix $B$ and when we simply write $\varrho$ and $\varrho_1,$ their values will be always those defined in \eqref{r} and \eqref{r1} respectively.

In the rest of the paper,  $f(n)\sim g(n)$ means $\lim_{n\rto}f(n)/g(n)=1,$ $f(n)=o(g(n))$ means $\lim_{n\rto}f(n)/g(n)=0$
%$f(n)\asymp g(n)$ means $\exists 0<C_1\le C_2<\infty$ such that $C_1\le f(n)/g(n)\le C_2$ for $n$ large enough,
and
 unless otherwise specified, $c$ and $C$  are some universal strictly positive constants, which may change from line to line. We always assume that empty product equals identity and empty sum equals $0.$ The convention $\sqrt{-1}=\mathrm i$ will also be adopted.

 Now we are ready to state the main result.

 \begin{theorem}\label{pnec}
   Suppose that condition (B1) and one of the conditions (B2)$_{ a},$ (B2)$_{ b}$ and (B2)$_{ c}$ hold.
   Assume further  that $|\varrho_1|<1,$  $bd\ne a\theta$ and $\kappa_i>0,\varkappa_i>0,i=1,2$ are proper constants. Then for $i\in\{1,2\}$
   \begin{align}\label{pns}
     &P(\nu>n|Z_0=\mb e_i)\sim \frac{\kappa_i}{\sum_{k=1}^{n+1} \varrho(M_1)^{-1}\cdots \varrho(M_{k-1})^{-1}} \text{ as }n\rto.
   \end{align}
   Furthermore, if $\varrho_1\ne\frac{1}{2}\z(a+b+1-\sqrt{(a+b+1)^2+4\frac{bd- a\theta}{\theta-b}}\y),$  then for $i\in\{1,2\}$
   \begin{align}\label{pnse}
     &P(\nu=n|Z_0=\mb e_i)\sim\frac{\varkappa_i\varrho(M_1)^{-1}\cdots \varrho(M_n)^{-1}}{\z(\sum_{k=1}^{n+1} \varrho(M_1)^{-1}\cdots \varrho(M_{k-1})^{-1}\y)^2}, \text{ as }n\rto;
   \end{align}
   otherwise, if  $\varrho_1=\frac{1}{2}\z(a+b+1-\sqrt{(a+b+1)^2+4\frac{bd- a\theta}{\theta-b}}\y),$ then for $i\in\{1,2\},$
    \begin{align}\label{pnso}
     &P(\nu=n|Z_0=\mb e_i)=o\z(\frac{\varrho(M_1)^{-1}\cdots \varrho(M_n)^{-1}}{\big(\sum_{k=1}^{n+1} \varrho(M_1)^{-1}\cdots \varrho(M_{k-1})^{-1}\big)^2}\y), \text{ as }n\rto.
   \end{align}
 \end{theorem}
 \begin{remark}
 By some elementary computation, one can show that if $bd>a\theta,$ then $\varrho_1=\frac{1}{2}\z(a+b+1-\sqrt{(a+b+1)^2+4\frac{bd- a\theta}{\theta-b}}\y)$ if and only if $\theta=b+1,$ see \cite{wy} for details.
 \end{remark}

 The proof of Theorem \ref{pnec} relies on the asymptotics of the product of nonhomogeneous matrices $A_k,k\ge 1,$ which have their own interests. See the two theorems below.
 \begin{theorem}\label{mpr}
 Suppose that condition (B1) and one of the conditions (B2)$_{ a},$ (B2)$_{ b}$ and (B2)$_{ c}$ hold. Then for $i,j\in\{1,2\}$ and $k\ge 1,$ there exists a number $c(k,i,j)\ne0$ such that
 \begin{align}\label{rml}
   \lim_{m\rto}\frac{\mathbf{e}_i A_{k} \cdots A_{m}  \mathbf e_j^t}{\varrho(A_{k}) \cdots \varrho(A_{m})}=c(k,i,j).
 \end{align}
 \end{theorem}
 \begin{remark}
   In \cite{hs20}, similarly result is proved if for each $k\ge1,$ all elements of $A_k$ are nonnegative, see Theorem 1 therein. But now, under condition (B1), $\tilde d_k$ might be negative. So, though the basic idea is somewhat the same, the proof of Theorem \ref{mpr} differs on a large extent  from the one of Theorem 1 in \cite{hs20}.
 \end{remark}
 \begin{theorem}\label{s}
  Suppose that condition (B1) and one of the conditions (B2)$_{ a},$ (B2)$_{ b}$ and (B2)$_{ c}$ hold.
   If  $bd\ne a\theta,$ then
   \begin{align}\label{sa}
     \sum_{k=1}^{n+1}\mb e_1A_k\cdots A_n\mb e_1^t \sim c\sum_{k=1}^{n+1}\varrho(A_k)\cdots\varrho(A_{n}), \text{ as }n\rto.
   \end{align}
\end{theorem}
\begin{remark}
  Although we have shown in Theorem \ref{mpr} that  for each $k\ge1,$ $\exists c(k)\ne0$ such that $\mb e_1A_k\cdots A_n\mb e_1^t\sim c(k) \varrho(A_k)\cdots\varrho(A_{n})$ as $n\rto,$  it is not an easy task to prove Theorem \ref{s}, because that the constants $c(k),k\ge1$ involved there are mutually different and moreover every summand in $\sum_{k=1}^{n+1}\mb e_1A_k\cdots A_n\mb e_1^t$ depends on $n.$
\end{remark}

\noindent{{\bf 1.3 Outline of the paper.} The remainder of the paper is organized as follows. In Section \ref{pr}, we give some preliminary results of continued fractions and also some facts on products of the matrices $M_k$  and $A_k.$ Section \ref{ae} and Section \ref{pc} are  devoted to proving Theorem \ref{mpr} and Theorem \ref{s} respectively. Finally, in Section \ref{pm}, based on Theorem \ref{mpr} and Theorem \ref{s}, we finish the proof  of Theorem \ref{pnec}.

\section{Preliminary results}\label{pr}
 Products of 2-by-2 matrices are closely related to continued fractions and therefore continued fractions are important tools to prove Theorem \ref{pnec}. To begin with, we introduce some basics of continued fractions.
 \subsection{Continued fractions and their tails}

Let $\beta_k,\alpha_k,k\ge 1$ be certain  real numbers. For $1\le k\le n,$ We denote by
\begin{equation}\label{aprx}
\xi_{k,n}\equiv\frac{\beta_k}{\alpha_k}\begin{array}{c}
                                \\
                               +
                             \end{array}\frac{\beta_{k+1}}{ \alpha_{k+1}}\begin{array}{c}
                                \\
                               +\cdots+
                             \end{array}\frac{\beta_n}{\alpha_n}:=\dfrac{\beta_k}{\alpha_k+\dfrac{\beta_{k+1}}{\alpha_{k+1}+_{\ddots_{\textstyle +\frac{\textstyle\beta_{n}}{\textstyle\alpha_{n} } }}}}
\end{equation}
 the $(n-k+1)$-th approximant of a continued fraction
 \begin{align}\label{xic}
   &\xi_k:=\frac{\beta_{k}}{\alpha_{k}}\begin{array}{c}
                                \\
                               +
                             \end{array}\frac{\beta_{k+1}}{\alpha_{k+1 }}\begin{array}{c}
                                \\
                               +
                             \end{array}\frac{\beta_{k+2}}{\alpha_{k+2}}\begin{array}{c}
                                \\
                               +\cdots
                             \end{array}.
\end{align}
We call $\xi_k,k\ge1$ in \eqref{xic}  the tails and $h_k:= \frac{\beta_k}{\alpha_{k-1}}_{+}\frac{\beta_{k-1}}{ \alpha_{k-2}}_{+\cdots+}\frac{\beta_2}{\alpha_1},k
                             \ge2$
the critical tails of the continued fraction ${\frac{\beta_{1}}{\alpha_{1}}}_{+}\frac{\beta_{2}}{\alpha_{2 }}_{+\cdots}$ respectively. We remark that in the literature, the $n$-th tail of a continued fraction  ${\frac{\beta_{1}}{\alpha_{1}}}_{+}\frac{\beta_{2}}{\alpha_{2 }}_{+\cdots}$ is usually denoted by $f^{(n)}={\frac{\beta_{n+1}}{\alpha_{n+1}}}_{+}\frac{\beta_{n+2}}{\alpha_{n+2 }}_{+\cdots}$
 and the critical tails are also slightly different from $h_n,n\ge1$ above.

If  $\lim_{n\rto}\xi_{k,n}$ exists, we say that the continued fraction $\xi_k$  is convergent and its value is defined as $\lim_{n\rto}\xi_{k,n}.$
The lemma below will be use times and again.
 \begin{lemma}\label{ct}
If $
  \lim_{n\rto}\alpha_n=\alpha\ne0, \lim_{n\rto}\beta_n=\beta, \text{ and } \alpha^2+4\beta\ge0,$
 then
 for any $k\ge1,$  $\lim_{n\rto}\xi_{k,n}$ exists and furthermore \begin{align}
  \lim_{k\rto}h_k= \lim_{k\rto}\xi_k=\frac{\alpha}{2}\z(\sqrt{1+4\beta/\alpha^2}-1\y).\no
 \end{align}
 \end{lemma}
The proof Lemma \ref{ct} can be found in many references,
  we refer the reader to \cite{lor}, see  discussion between (4.1) and (4.2) on page 81 therein.

\subsection{Some facts on  matrices $A_k$ and $M_k$}

By assumptions on the numbers $a_k,b_k,d_k$ and $\theta_k,$ we have
\begin{align}
  A_k\cdots A_n&=\Lambda_k^{-1}M_k\cdots M_n\Lambda_{n+1},n\ge k\ge1,\no\\
  \mb e_1 A_k\cdots A_n \mb e_1^t&=\mb e_1M_k\cdots M_n(1,\theta_{n+1}/b_{n+1})^t>0,n>k\ge1,\label{apm1}\\
  \mb e_1 A_k\cdots A_n \mb e_2^t&=\mb e_1M_k\cdots M_n\mb e_2^t>0,n> k\ge1,\label{apm2}\\
  \mb e_2 A_k\cdots A_n \mb e_1^t&=(-{\theta_k}/{b_k},1)M_k\cdots M_n(1,\theta_{n+1}/b_{n+1})^t,n\ge k\ge1. \label{apm21}
    \end{align}
    Note that the conditions of Theorem \ref{pnec} can not ensure  $\Delta_k:=\Big(a_k+\frac{b_k\theta_{k+1}}{b_{k+1}}\Big)^2+4(b_kd_k-a_k\theta_k)>0$ for each $k\ge1.$ Therefore, the matrix $A_k$ may possess complex eigenvalues. However, since $\lim_{k\rto}\Delta_k=(a-\theta)^2+4bd>0,$ there exists a number $N_1>0$ such that   $\Delta_k>0,$ $\forall k\ge N_1,$ and thus
  \begin{align}\label{rhk}
         \varrho(A_k)&=\frac{\sqrt{\tilde a_k^2+4\tilde b_k\tilde d_k}+\tilde a_k}{2}
         \\&=\frac{1}{2}\z(a_k+\frac{b_k\theta_{k+1}}{b_{k+1}}+\sqrt{\Big(a_k+\frac{b_k\theta_{k+1}}{b_{k+1}}\Big)^2+4(b_kd_k-a_k\theta_k)}\y),\forall k\ge N_1.\no
    %&=\frac{1}{2}\Big(a_k+\theta_k+b_k\Delta_k+\sqrt{\big(a_k+\theta_k+b_k\Delta_k\big)^2+4(b_kd_k-a_k\theta_k)}\Big)>0.\no
  \end{align}
 Furthermore, if $bd\ne a\theta,$ then $\lim_{n\rto}\tilde d_n=d-a\theta/b\ne0.$ Thus, if $bd<a\theta,$ then
 \begin{equation*}\begin{split}
   &\exists \varepsilon>0 \text{ and }N_2>0, \text{ such that }  \tilde d_n<-\varepsilon, \ \forall n\ge N_2 \text{  and consequently }\\
& \mb e_2 A_k\cdots A_n \mb e_1^t<0,n\ge k\ge N_2 \text{ and }\mb e_2 A_k\cdots A_n \mb e_2^t<0,n>k\ge N_2;\end{split}
\end{equation*}
otherwise, if $bd>a\theta,$ then
 \begin{equation*}\begin{split}
   &\exists \varepsilon>0 \text{ and }N_2>0, \text{ such that }  \tilde d_n>\varepsilon, \ \forall n\ge N_2 \text{  and consequently }\\
&
 \mb e_2 A_k\cdots A_n \mb e_1^t>0,n\ge k\ge N_2 \text{ and }\mb e_2 A_k\cdots A_n \mb e_2^t>0,n>k\ge N_2.\end{split}
\end{equation*}
\section{Asymptotics of elements of $\prod_{j=k}^{m}A_j$}\label{ae}
The main task of this section is to prove Theorem \ref{mpr}. The proof is based on the spectral radius estimation derived in \cite{235} and some analysis of the fluctuations of the tails and critical tails of continued fractions.

\subsection{Lower and upper bounds for $\frac{\mathbf{e}_1 A_{k} \cdots A_{m}  \mathbf e_1^t}{\varrho(A_{k}) \cdots \varrho(A_{m})}$ }
%In this subsection, we show that for $k$ large enough, $\frac{\mathbf{e}_1 A_{k} \cdots A_{m}  \mathbf e_1^t}{\varrho(A_{k}) \cdots \varrho(A_{m})},m\ge k$ are uniformly bounded away from $0$ and infinity.
The main result of this subsection is the following lemma.
\begin{lemma}\label{lwb}
  Suppose that  condition (B1) holds. Then  $\exists N_0>0$ such that for each $k\ge N_0,$ $C_k^{-1}<\frac{\mathbf{e}_1 A_{k} \cdots A_{m}  \mathbf e_1^t}{\varrho(A_{k}) \cdots \varrho(A_{m})}<C_k, \forall m\ge k,$ where $C_k>0$ is a proper number.
\end{lemma}
Lemma \ref{lwb} is a direct consequence of  the following four auxiliary  lemmas.

\begin{lemma}\label{bma}
 Suppose the condition (B1) is satisfied.  Then with $N_1$ the one in \eqref{rhk},
 for each $k\ge N_1,$ $\exists$ $ \mathcal C_k>0$ and $\mathcal B_k>0$ such that \begin{align*}
   &\prod_{j=k}^{n}\varrho(A_j)\sim \mathcal C_k\prod_{j=k}^{n}\varrho(M_j)\text{ and } \sum_{j=k}^{n+1}\prod_{i=k}^{j-1} \varrho(A_i)^{-1}\sim \mathcal B_k\sum_{j=k}^{n+1}\prod_{i=k}^{j-1} \varrho(M_i)^{-1},\text{ as }n\rto.\end{align*}
  \end{lemma}
  For the proof of Lemma \ref{bma}, we refer the reader to \cite{wy}.

   \begin{lemma}\label{rfm}
    Suppose that  condition (B1) holds.  Then there exist constants $0<\zeta<\gamma<\infty$ such that for $m\ge k\ge1,$
    $
 \zeta \leq \frac {\varrho(M_{k} \cdots M_{m})}{\varrho(M_{k}) \cdots \varrho(M_{m})}\leq \gamma.
$

  \end{lemma}

\proof For vectors $\mathbf v= \left(
        \begin{array}{c}
        v_1\\
        v_2  \\
        \end{array}
       \right)$ and $\mathbf u=\left(
        \begin{array}{c}
        u_1\\
        u_2  \\
        \end{array}
       \right),$ set $$\frac{\mathbf v}{\mathbf u}:= \left(
        \begin{array}{c}
        v_1/u_1\\
        v_2/u_2  \\
        \end{array}
       \right),\ \left(\frac{\mathbf v}{\mathbf u}\right)_{\textrm{min}}:=\min\z\{\frac{v_1}{u_1},\frac{v_2}{u_2}\y\}, \ \left(\frac{\mathbf v}{\mathbf u}\right)_{\textrm{max}}:=\max\z\{\frac{v_1}{u_1},\frac{v_2}{u_2}\y\}.$$
  Let $\mathbf v_n$ be a right eigenvector of $A_n$ corresponding to the largest eigenvalue $\varrho(A_n)$. Then we can choose $\mathbf v_n$ to be
$
\mathbf{v}_n=\left( \varrho(A_n)-\theta_n,d_n\right)^t.
$
For $ m\ge k\ge 1,$ write \begin{align*}
  &\gamma_{k,m}:=\Big(\frac{\mathbf v_k}{\mathbf v_{k-1}}\Big)_{\max} \cdots \Big(\frac{\mathbf v_{m+1}}{\mathbf v_{m}}\Big)_{\max}
\Big(\frac{\mathbf v_{m}}{\mathbf v_{k}}\Big)_{\max},\\
&\zeta_{k,m}:=\Big(\frac{\mathbf v_k}{\mathbf v_{k-1}}\Big)_{\min} \cdots \Big(\frac{\mathbf v_{m+1}}{\mathbf v_{m}}\Big)_{\min}
\Big(\frac{\mathbf v_{m}}{\mathbf v_{k}}\Big)_{\min}.
\end{align*}
Applying \cite[Theorem 1, page 228]{235}, for $k\ge1,$ we have
\begin{align*}
\zeta_{k,m}\leq
\frac {\varrho(M_{k} \cdots M_{m})}{\varrho(M_{k}) \cdots \varrho(M_{m})}\le \gamma_{k,m}.
\end{align*}
It remains to show that both $\zeta_{k,m}^{-1}$ and $\gamma_{k,m},m\ge k\ge1$ are uniformly bounded away from $\infty.$  To this end, set $\epsilon_n=(\mathbf v_n/\mathbf v_{n-1})_{\max}-1, n\ge2.$ Then by condition (B1), for $n\ge 2,$ we have
\begin{align*}
  |\epsilon_n|&\le\max\left\{\left|\frac{\varrho(M_n)-\theta_n}{\varrho(M_{n-1})-\theta_{n-1}}-1\right|, \left|\frac{d_n}{d_{n-1}}-1\right|\right\}\\
  &  \le \left|\frac{\varrho(M_n)-\theta_n}{\varrho(M_{n-1})-\theta_{n-1}}-1\right|+ \left|\frac{d_n}{d_{n-1}}-1\right|\\
  &\le c(|\varrho(M_{n})-\varrho(M_{n-1})|+|d_{n}-d_{n-1}|+|\theta_n-\theta_{n-1}|)\\
  &\le c(|a_n-a_{n-1}|+|b_n-b_{n-1}|+|d_n-d_{n-1}|+|\theta_{n}-\theta_{n-1}|) <\infty ,
\end{align*} which implies $\sum_{n=2}^\infty \log(1+|\epsilon_n|)<\infty.$ As a consequence,  $$\gamma_{k,m}\le \max\Big\{\frac{\varrho(M_m)-\theta_m}{\varrho(M_k)-\theta_k},
\frac{d_m}{d}\Big\}\prod_{n=2}^\infty(1+|\epsilon_n|)<\gamma$$ for some number $\gamma<\infty$ independent of $k$ and $m.$

Since $\zeta_{k,m}^{-1}=\Big(\frac{\mathbf v_{k-1}}{\mathbf v_k}\Big)_{\max} \cdots \Big(\frac{\mathbf v_{m}}{\mathbf v_{m+1}}\Big)_{\max}
\Big(\frac{\mathbf v_{k}}{\mathbf v_{m}}\Big)_{\max},$ a similar argument also yields that
$\zeta_{k,m}^{-1},m\ge k\ge1$ are uniformly bounded away from $\infty.$ Consequently, Lemma \ref{rfm} is proved.\qed
%\begin{remark}\label{brb}
%  If  condition (B1) is replaced by ``(B1)$'$: for some numbers $a,b,d>0,$  $a_k\rightarrow a,$ $b_k\rightarrow b,$  $d_k\rightarrow d$  as $k\rto,$ and in addition, $\frac{a_k}{d_k}, \frac{b_k}{d_k}$ are increasing (or decreasing) in $k$ simultaneously", then we may choose a right eigenvector corresponding to  $\varrho(A_n)$ as
%$
%\mathbf{v}_n=\left(
%        \varrho(A_n)/d_n,
%        1  \right)^t.
%$ Let $\gamma_k$ and $\zeta_k$ be those defined in \eqref{ga} and \eqref{ze}. Under condition (B1)$'$, $\varrho(A_n)/d_n$ is monotone in $n.$ If it is increasing in $n,$ then
%$\gamma_k=\frac{\varrho(A_k)d_1}{\varrho(A_1)d_k}$ and $\zeta_k=\frac{\varrho(A_1)d_k}{\varrho(A_k)d_1},k\ge 1.$
%Since $\lim_{k\rto}\gamma_k=\lim_{k\rto}\zeta_k^{-1}=c$ for some $0<c<\infty,$ then both $\gamma_k,k\ge1$ and $\zeta_k,k\ge1$ are uniformly bounded away from $0$ and infinity.
%Otherwise, if $\varrho(A_n)/d_n$ is decreasing  in $n,$  things can be done by a similar approach.
%\end{remark}

  \begin{lemma} \label{ral} Suppose that  condition (B1) holds and for $k\ge1,$ let $\xi_k$ be the one in \eqref{xic} with $\alpha_n=\frac{\tilde a_n}{\tilde b_n},\beta_n =\frac{\tilde d_n}{\tilde b_n},n\ge1.$ Then $\exists N_3>0$ such that
\begin{align}
%&\lim_{k\rto}\xi_k=-\frac{\varrho_1}{b},\label{lf}\\
&\lim_{m\rto}\frac{\varrho(A_{k} \cdots A_{m})}{\mathbf{e}_1 A_{k} \cdots A_{m}  \mathbf e_1^t}=1+b\varrho^{-1}\xi_k>c>0, k\ge N_3. \label{kmp}
 \end{align}
  \end{lemma}

\proof
  For $m\ge k\ge 1,$ write
  $A_{k,m}:=A_k\cdots A_m=\z( \begin{array}{cc}
                        A_{k,m}(11) & A_{k,m}(12) \\
                        A_{k,m}(21)& A_{k,m}(22)
                      \end{array}
  \y).$
%
%
%Indeed, an application of the ergodicity theorem of the product of nonnegative matrices(see \cite[Theorem 3.3]{se}) yields the existence of the limits  and the second equality in both (\ref{fm}) and (\ref{gm}).
% To compute $ f$ and $ g_m,$
Noticing that
\begin{align*}
A_{k,m}&=\left(
      \begin{array}{cc}
      \tilde a_{k} & \tilde b_{k} \\
      \tilde d_k & 0 \\
      \end{array}
    \right) \cdots \left(
      \begin{array}{cc}
      \tilde a_{m} & \tilde b_{m} \\
      \tilde d_m & 0 \\
      \end{array}
    \right)=\tilde b_k \cdots \tilde b_{m} \left (
    \begin{array}{cc}
      \alpha_k & 1 \\
     \beta_k & 0 \\
      \end{array}
    \right) \cdots \left(
      \begin{array}{cc}
      \alpha_m & 1 \\
      \beta_m & 0 \\
      \end{array}
    \right),
    \end{align*}
    thus we get
    \begin{align*}
      \frac{A_{k,m+1}(22)}{A_{k,m+1}(12)}=\frac{A_{k,m}(21)}{A_{k,m}(11)},
    \end{align*}
and by forward and backward induction we have
\begin{align}
&\frac{A_{k,m}(21)}{A_{k,m}(11)}=\frac{\beta_{k}}{\alpha_k}\underset{+}\quad \frac{\beta_{k+1}}{\alpha_{k+1}}\underset{+}\quad \underset{\cdots+}\quad \frac{\beta_{m-1}}{\alpha_{m-1}}\underset{+}\quad \frac{\beta_{m}}{\alpha_m}=:\xi_{k,m}\label{cfaa},\\
&\frac{A_{k,m}(12)}{A_{k,m}(11)}=\frac{1}{\alpha_m}\underset{+}\quad \frac{\beta_{m}}{\alpha_{m-1}}\underset{+}\quad \frac{\beta_{m-1}}{\alpha_{m-2}}\underset{+}\quad \underset{\cdots+}\quad \frac{\beta_{k+2}}{\alpha_{k+1}}\underset{+}\quad \frac{\beta_{k+1}}{\alpha_k}=:h_{m,k},\label{cfb}\\
&\frac{A_{k,m}(22)}{A_{k,m}(21)}=\frac{1}{\alpha_m}\underset{+}\quad \frac{\beta_{m}}{\alpha_{m-1}}\underset{+}\quad \frac{\beta_{m-1}}{\alpha_{m-2}}\underset{+}\quad \underset{\cdots+}\quad \frac{\beta_{k+3}}{\alpha_{k+2}}\underset{+}\quad \frac{\beta_{k+2}}{\alpha_{k+1}}.\no
\end{align}
Since $\alpha_k\rightarrow(a+\theta)/b$ and  $\beta_k\rightarrow(bd-a\theta)/b^2$ as $n\rto,$
applying  Lemma \ref{ct} we get
 \begin{align}
  &\lim_{m\rto}\frac{A_{k,m}(21)}{A_{k,m}(11)}=\lim_{k\rto}\frac{A_{k,m}(22)}{A_{k,m}(12)}=\xi_k,\label{cfa}\\
  &\lim_{m\rto}\frac{A_{k,m}(12)}{A_{k,m}(11)}=\lim_{k\rto}\frac{A_{k,m}(22)}{A_{k,m}(21)}
  %&=\frac{2b}{a+\theta+\sqrt{(a+\theta)^2+4(bd-a\theta)}}
  =b\varrho^{-1},\label{gm}\\
&\lim_{m\rto}\xi_k=-\frac{\varrho_1}{b}.\label{lxka}
\end{align}
We are ready to prove \eqref{kmp}. Write
\begin{align}\label{pkn}
  P_{k,m}&=A_{k,m}(12)A_{k,m}(21)-A_{k,m}(11)A_{k,m}(22).
\end{align}
If $(A_{k,m}(11)+A_{k,m}(22))^2+4P_{k,m}\ge 0,$ then
\begin{align}\label{ra}
  &\varrho(A_{k,m})=\frac{A_{k,m}(11)+A_{k,m}(22)}{2}+\frac{\sqrt{(A_{k,m}(11)+A_{k,m}(22))^2+4P_{k,m}}}{2};
\end{align}
otherwise, if
If $(A_{k,m}(11)+A_{k,m}(22))^2+4P_{k,m}< 0,$ then
\begin{align*}
  &\varrho(A_{k,m})=\frac{1}{2}\sqrt{(A_{k,m}(11)+A_{k,m}(22))^2+|(A_{k,m}(11)+A_{k,m}(22))^2+4P_{k,m}|}.
\end{align*}
Note that $\forall m\ge k\ge1,$ $A_{k,m}(11)>0$ and by \eqref{cfa}-\eqref{lxka}, we have
\begin{align}
&\lim_{k\rto}(1+b\varrho^{-1}\xi_k)=1-\frac{\varrho_1}{\varrho}>0,\label{lk}\\
  &\lim_{m\rto}\frac{(A_{k,m}(11)+A_{k,m}(22))^2+4P_{k,m}}{(A_{k,m}(11))^2}=\z(1+b\varrho^{-1}\xi_k\y)^2.\no
\end{align}
Therefore there exists $N_3>0,$ such that  for each $k>N_3,$  we have $1+b\varrho^{-1}\xi_k>0,$ and for $m$ large enough
the eigenvalues of $A_{k,m}$ are real so that \eqref{ra} holds. Thus,
  taking \eqref{cfa}-\eqref{lxka} into account, for $k\ge N_3,$  we get \begin{align*}
\lim_{m\rto}&\frac{\varrho(A_{k} \cdots A_{m})}{\mathbf{e}_1 A_{k} \cdots A_{m}  \mathbf e_1^t}=\lim_{m \rto}\frac{\varrho(A_{k,m})}{A_{k,m}(11)}\no\\
&=1/2(1+b\varrho^{-1}\xi_k+|1+b\varrho^{-1}\xi_k|)=1+b\varrho^{-1}\xi_k>c>0,
%\lim_{m\rto}&\frac{\varrho(A_{k} \cdots A_{m})}{\mathbf{e}_2 A_{k} \cdots A_{m}  \mathbf e_1^t}=\lim_{m \rto}\frac{\varrho(A_{k,m})}{A_{k,m}(21)}\no\\
%&=1/2 \xi_k^{-1}(1+b\varrho^{-1}\xi_k+|1+b\varrho^{-1}\xi_k|)=\xi_k^{-1}(1+b\varrho^{-1}\xi_k)\ne0.
 \end{align*} which proves \eqref{kmp}.
 %Similarly, for general $k\ge1,$  dividing \eqref{ra} and \eqref{rab} by $A_{k,m}(11)$ and letting $m\rto,$ we get \eqref{ll} and \eqref{ul}.
 %\begin{align*}
%\frac{\sqrt{2}}{2}(1+b\varrho^{-1}\xi_k)\le \varliminf_{n\rto}\frac{\varrho(A_{k,m})}{A_{k,m}(11)}\le \varlimsup_{n\rto}\frac{\varrho(A_{k,m})}{A_{k,m}(11)}=1+b\varrho^{-1}\xi_k.
% \end{align*}
 The lemma is proved. \qed
\begin{lemma}\label{ff} For $k\ge1,$ set $\varphi_k=:\frac{1+b\varrho^{-1}\xi_k}{1+b\varrho^{-1}\xi_k+b\varrho^{-1}\big(\frac{\theta_k}{b_k}-\frac{\theta}{b}\big)}.$ Then, $\exists N_4>0$ such that   \begin{align}\label{rma}\lim_{n\rto}\frac{\varrho(A_k\cdots A_n)}{\varrho(M_k\cdots M_n)}=\varphi_k>c>0,k\ge N_4.\end{align}
\end{lemma}
\proof Write $A_{k,n}=A_k\cdots A_n,$ let $P_{k,n}$ be the one in \eqref{pkn} and set \begin{align*}
  Q_{k,n}&=A_{k,n}(11)+A_{k,n}(22)+A_{k,n}(12)\Big(\frac{\theta_k}{b_k}-\frac{\theta_{n+1}}{b_{n+1}}\Big).
  \end{align*}
Since $M_k\cdots M_n =\Lambda_kA_{k,n}\Lambda_{n+1}^{-1},$ then
\begin{align*}
  &\varrho(M_k\cdots M_n)=\frac{1}{2}\z(Q_{k,n}+\sqrt{Q_{k,n}^2+4P_{k,n}}\y).
\end{align*}
Using \eqref{cfa} and \eqref{gm}, we get
\begin{align}\label{rmb}
  \lim_{n\rto}\frac{\varrho(M_k\cdots M_n)}{A_{k,n}(11)}=1+b\varrho^{-1}\xi_k+b\varrho^{-1}\Big(\frac{\theta_k}{b_k}-\frac{\theta}{b}\Big).
\end{align}
%Since $\lim_{k\rto}(1+b\varrho^{-1}\xi_k)=1-\frac{\varrho_1}{\varrho}>0$ and $\lim_{k\rto}\theta_k/b_k\rightarrow \theta/b,$ then $$\lim_{k\rto}\frac{1/2(1+b\varrho^{-1}\xi_k+|1+b\varrho^{-1}\xi_k|)}{1+b\varrho^{-1}\xi_k+b\varrho^{-1}\big(\frac{\theta_k}{b_k}-\frac{\theta}{b}\big)}=1.$$
Let $N_3$ be the one in Lemma \ref{ral}. Using \eqref{lk} and the fact $\lim_{k\rto}\theta_k/b_k\rightarrow \theta/b,$ we have $\lim_{k\rto}\varphi_k=1.$ Hence, $\exists N_4>N_3$ such that for $k>N_4,$ $\varphi_k>1/2.$
As a result, \eqref{rmb} together with \eqref{kmp} implies that \eqref{rma} is true.
  \qed

\subsection{Fluctuations of tails and critical tails of continued fractions}
For $k\ge1,$ let \begin{align}\label{fk}
  f_k&=\frac{ \tilde b_k \tilde d_k^{-1}}{ \tilde a_k\tilde d_k^{-1}}\begin{array}{c}
                                \\
                               +
                             \end{array}\frac{ \tilde b_{k-1} \tilde d_{k-1}^{-1}}{ \tilde a_{k-1}\tilde d_{k-1}^{-1}}\begin{array}{c}
                                \\
                               +\cdots+
                             \end{array}\frac{ \tilde b_1 \tilde d_1^{-1}}{ \tilde a_1 \tilde d_1^{-1}},\\
\xi_k&=\frac{ \tilde b_k^{-1} \tilde d_{k+1}^{-1}}{ \tilde a_k \tilde b_k^{-1} \tilde d_{k+1}^{-1}}\begin{array}{c}
                                \\
                               +
                             \end{array}\frac{ \tilde b_{k+1}^{-1} \tilde d_{k+2}^{-1}}{ \tilde a_{k+1}\tilde b_{k+1}^{-1}\tilde d_{k+2}^{-1}}\begin{array}{c}
                                \\
                               +\cdots.
                             \end{array}\no
\end{align}
Set also \begin{align}
  \varepsilon^f_k&=f_k-\tilde b_{k+1}\varrho(A_{k+1})^{-1},\ \varepsilon^\xi_k=\xi_k-\varrho(A_{k})^{-1}, k\ge1,\label{de}\\
  \delta_k^f&=\tilde b_kd_k^{-1}-\tilde b_{k+1}\varrho(A_{k+1})^{-1}(\tilde a_kd_k^{-1}+\tilde b_{k}\varrho(A_{k})^{-1}),k\ge2,\no\\
  \delta_k^\xi&=\tilde b_k^{-1}\tilde d_{k+1}^{-1}-\varrho(A_k)^{-1}(\tilde a_k\tilde b_k^{-1}\tilde d_{k+1}^{-1}+\varrho(A_{k+1})^{-1}),k\ge1.\no
\end{align}
Suppose that condition (B1) holds. Then applying Lemma \ref{ct}, we have
\begin{align}%\label{xkl}
  f_k\rightarrow b\varrho^{-1},\ \xi_k\rightarrow\varrho^{-1} \text{ and hence }\varepsilon_k^f \rightarrow0,\ \varepsilon_k^\xi \rightarrow0 \text{ as }k\rto.\no
\end{align}

The following lemma, which can be proved by some arguments similar to the proofs of \cite[Lemma 4]{hs20}  and \cite[Lemma 12]{wy}, studies the fluctuations of both $\varepsilon_k^f$ and $\varepsilon_k^\xi,$ $k\ge1.$
\begin{lemma}\label{dxf}Suppose that condition (B1) and one of the conditions (B2)$_{a},$ (B2)$_{ b}$  and  (B2)$_{c}$ hold. Then there exists some number $q$ with $|q|\le1$ such that
\begin{align*}%\label{qdd}
  &\lim_{k\rto}{\delta_{k+1}^f}/{\delta_k^f}=\lim_{k\rto}{\delta_{k+1}^\xi}/{\delta_k^\xi}=q,\\
 %\label{qff}
 &\lim_{k\rto}\frac{\varepsilon_{k+1}^\xi}{\varepsilon_{k}^\xi}=q,\  \lim_{k\rto}\frac{\varepsilon_{k+1}^f}{\varepsilon_{k}^f}=q\text{ or }\frac{a\theta- bd}{\varrho^2}.
\end{align*}
\end{lemma}

\subsection{Proof of Theorem \ref{mpr}.}
 For $m\ge k\ge1$  we write
$x_{k,m}:=\frac{\mathbf e_1A_k\cdots A_m\mathbf e_1^t}{\varrho(A_k)\cdots\varrho(A_m)}$ for simplicity. In view of (\ref{cfa}) and (\ref{gm}), to prove \eqref{rml}, it suffices to show that $x_{k,m}\rightarrow c(k)$ as $m\rto$ for some $0<c(k)<\infty.$
%
%It follows from Lemma \ref{lwb} that there exist some constants $0<c_3<c_4<\infty$ independent of $k$ and $m$ such that
%$ c_3\le x_{k,m}\le c_4,\forall k\ge m\ge1.$
%Therefore, if \eqref{rml} is true,  we must have  $c_3\le c(k)\le c_4,k\ge1.$ Thus the second part of Theorem \ref{mpr} holds.

%{\it Step 1: We show that $x_{k,m}\rightarrow c(k) $ as $m\rto$ for some $0<c(k)<\infty.$}

To begin with, we show first that for fixed $j\ge1,$ $x_{j,m},m\ge j$ are uniformly bounded away from $0$ and infinity. Indeed, using again the notation $A_{m,n}=A_m\cdots A_n,$
 in view of \eqref{cfaa} and \eqref{cfb}, for $m\ge n\ge j\ge1,$
\begin{align}\no
  x_{j,m}&=\frac{\mb e_1A_{j}\cdots A_m\mb e_1^t}{\varrho(A_j)\cdots \varrho(A_m)}\\
 & =\frac{A_{j,n}(11)}{\prod_{i=j}^{n}\varrho(A_i)}\frac{A_{n+1,m}(11)}{\prod_{i=n+1}^m\varrho(A_{i}) }+
  \frac{A_{j,n}(12)}{\prod_{i=j}^{n}\varrho(A_i)}\frac{A_{n+1,m}(21)}{\prod_{i=n+1}^m\varrho(A_{i})}\no\\
  &=\frac{A_{j,n}(11)}{\prod_{i=j}^{n}\varrho(A_i)}\frac{A_{n+1,m}(11)}{\prod_{i=n+1}^m\varrho(A_{i}) }\z(1+\frac{A_{j,n}(12)}{A_{j,n}(11)}\frac{A_{n+1,m}(21)}{A_{n+1,m}(11)}\y)\no\\
  &=x_{j,n}x_{n+1,m}(1+h_{n,j}\xi_{n+1,m}).\no
\end{align}
Let $\xi_i,i\ge1$ be those in Lemma \ref{ral}. Then by \eqref{cfa}-\eqref{lxka} we get
\begin{align*}
  &\lim_{m\rto}(1+h_{n,j}\xi_{n+1,m})=1+h_{n,j}\xi_{n+1},\\
   &\lim_{n\rto}(1+h_{n,j}\xi_{n+1})=1-\frac{\varrho_1}{\varrho}\in(0,2).
\end{align*}
Thus, we can find a number $N>N_0$ such that $0<c<1+h_{n,j}\xi_{n+1}<2-c, \forall n\ge N,$ where  $N_0$ is the one in Lemma \ref{lwb}.

  Now fix $k>N.$ There is a number $N_5>k$ such that for $m>N_5,$
 $$C^{-1}<1/2 (1+h_{k,j}\xi_{k+1})<1+h_{k,j}\xi_{k+1,m}<2(1+h_{k,j}\xi_{k+1})<C.$$
 Then \begin{align}\label{xjm}
   C^{-1}x_{j,k}x_{k+1,m}<x_{j,m}<Cx_{j,k}x_{k+1,m},\forall m>N_5,j\ge1.
 \end{align}
 But by Lemma \ref{lwb}, we have $C_{k+1}^{-1}<x_{k+1,m}<C_{k+1},\forall m\ge k+1,$ with $C_{k+1}>0$ certain constant.
 Let $\zeta_{k+1}=\max\{x_{k+1,m}, k+1\le m\le N_5\}$ and set $C_{j,k}=(Cx_{j,k}(\zeta_{k+1}\vee C_{k+1}))\vee \max_{j\le m\le N_5}x_{j,m},$ where $C$ is the constant in \eqref{xjm}.
 Then we have \begin{align}\label{xlu}
   C_{j,k}^{-1}<x_{j,m}<C_{j,k}, \forall m\ge j\ge1.
 \end{align}

 With \eqref{xlu} and Lemma \ref{dxf} in hands, the rest of the proof of Theorem \ref{mpr} is more or less similar to the one of Theorem 1 in \cite{hs20}, but there are some differences in details. For convenience of the reader,  we provide a complete proof here.

Next, we prove only $x_{1,m}\rightarrow c $ as $m\rto$ for some $0<c<\infty,$  since for $k\ge2$ the convergence of $x_{k,m}$ as $m\rto$ can be proved similarly. For simplicity, we write $x_{1,m}$ as $x_m$ and taking \eqref{xlu} into account, we keep always in mind that
\begin{equation}\label{xul}
  x_{m}=\frac{\mathbf e_1A_1\cdots A_m\mathbf e_1^t}{\varrho(A_1)\cdots\varrho(A_m)} \text{ and } \forall m\ge1, c_3\le x_m\le c_4,
\end{equation}
where $0< c_3\le c_4<\infty$ are  proper constants.

For $n\ge1,$ writing $f_n=\frac{\mb e_1A_1\cdots A_{n}\mathbf e_2^t}{\mb e_1A_{1}\cdots A_n\mathbf e_1^t},$ then
\begin{align}\label{fpa}
 f_n=\frac{\tilde b_n\mb e_1A_1\cdots A_{n-1}\mb e_1^t}{\mb e_1A_{1}\cdots A_{n-1}(\tilde a_n\mathbf e_1^t+\tilde d_n\mb e_2^t)}=\frac{\tilde b_n}{\tilde a_n+\tilde d_nf_{n-1}}.
\end{align}
Iterating \eqref{fpa}, we find that
$f_n,n\ge1$ coincide with the ones defined in (\ref{fk}).  Therefore, letting $\varepsilon_n^f,n\ge 1$ be the ones in \eqref{de},
 an application of  Lemma \ref{dxf} yields that for some number $q$ with $|q|\le1,$
\begin{align}\label{elk}
  \lim_{m\rto}\frac{\varepsilon_{m+1}^f}{\varepsilon_{m}^f}=q\text{ or }\frac{a\theta-bd}{\varrho^2}.
\end{align}
It is easy to check that $0<\frac{|a\theta-bd|}{\varrho^2}<1$ since $a\theta-bd\ne0.$

{\it Case 1: Suppose $|q|<1.$} Then $\lambda_0:=\max\{|q|,\frac{|a\theta-bd|}{\varrho^2}\}<1.$ Fix $0<\lambda<\lambda_0.$ By (\ref{elk}), there exists some $m_0>0$ such that $\z|\frac{\varepsilon_{m+1}^f}{\varepsilon_m^f}\y|\le \lambda $ for all $m>m_0.$ Hence
\begin{equation}\label{foa}
  \sum_{m=2}^{\infty}|\varepsilon_m^f|=\sum_{m=2}^\infty|f_m-\tilde b_{m+1}\varrho(A_{m+1})^{-1}| <\infty.
\end{equation}
Taking \eqref{rhk} into account, for $m\ge 1$ we have
\begin{align}\label{xf}
  x_{m+1}-x_m&=\frac{\tilde a_{m+1}-\varrho(A_{m+1})+\tilde d_{m+1}\frac{\mathbf e_1A_1\cdots A_m\mathbf e_2^t}{\mathbf e_1A_1\cdots A_m\mathbf e_1^t}}{\frac{\varrho(A_{1})\cdots\varrho(A_{m+1})}{\mathbf e_1A_1\cdots A_m\mathbf e_1^t}}\\
  &=\varrho(A_{m+1})^{-1}x_m(\tilde a_{m+1}-\varrho(A_{m+1})+\tilde d_{m+1}f_m)\no\\
  &=\varrho(A_{m+1})^{-1}x_m\tilde d_{m+1}(f_m-\tilde b_{m+1}\varrho(A_{m+1})^{-1})\no\\
  &=\varrho(A_{m+1})^{-1}x_m\tilde d_{m+1}\varepsilon_m^f.\no
\end{align}
Since $\varrho(A_{m}),m\ge 1$ are uniformly bounded away from $0$ and $\infty,$ then by (\ref{xul}), we have for some constant $0<c_5<\infty,$
\begin{align}\label{xda}
  |x_{m+1}-x_m|\le c_5 |\varepsilon_m^f|, \forall m\ge1.
\end{align}
Taking (\ref{xul}), (\ref{foa}) and (\ref{xda}) together, we conclude that for some constant $0<c<\infty,$ $\lim_{m\rto} x_m=c.$

{\it Case 2: Suppose $q=1$ and $ \lim_{m\rto}\frac{\varepsilon_{m+1}^f}{\varepsilon_{m}^f}=\frac{a\theta-bd}{\varrho^2}.$ } Since $\frac{|a\theta-bd|}{\varrho^2}<1,$ the proof goes exactly the same as Case 1.

{\it Case 3: Suppose $q=1$ and $\lim_{m\rto}\frac{\varepsilon_{m+1}^f}{\varepsilon_{m}^f}=q.$} Then there exists some number $m_1>0$ such that $\varepsilon_m^f=f_m-\tilde b_{m+1}\varrho(A_{m+1})^{-1}, m\ge m_1$ are all strictly positive or strictly negative, and consequently
\begin{align}\label{skof}\frac{\tilde a_{m+1}+\tilde d_{m+1}f_m}{\tilde a_{m+1}+\tilde d_{m+1}\tilde b_{m+1}\varrho(A_{m+1})^{-1}}<1(\text{or }>1), \text{ for all }m\ge m_1.\end{align}
But \begin{align}\label{ratio}
  \frac{x_{m+1}}{x_m}&=\frac{1}{\varrho(A_{m+1})}\z(\tilde a_{m+1}+\tilde d_{m+1}\frac{\mathbf e_1A_1\cdots A_m\mathbf e_2^t}{\mathbf e_1A_1\cdots A_m\mathbf e_1^t}\y)\\
  &=\frac{1}{\varrho(A_{m+1})}\z(\tilde a_{m+1}+\tilde d_{m+1}f_m\y)=\frac{\tilde a_{m+1}+\tilde d_{m+1}f_m}{\tilde a_{m+1}+\tilde d_{m+1}\tilde b_{m+1}\varrho(A_{m+1})^{-1}}.\no
\end{align}
Thus, by (\ref{skof}), $\frac{x_{m+1}}{x_m}<1(\text{or }>1)$ for all $m\ge m_1,$ that is, $x_m,m\ge m_1$ is monotone in $m.$ As a consequence, it follows from (\ref{xul}) that for some constant $0<c<\infty$ $\lim_{m\rto}x_m=c.$

{\it Case 4. Suppose that $q=-1$ and $ \lim_{m\rto}\frac{\varepsilon_{m+1}^f}{\varepsilon_{m}^f}=\frac{a\theta-bd}{\varrho^2}.$ } In this case, the proof is the same as Case 2.

{\it Case 5. Suppose that $q=-1$ and $ \lim_{m\rto}\frac{\varepsilon_{m+1}^f}{\varepsilon_{m}^f}=-1.$}
Combining (\ref{xf}) with \eqref{ratio}, we have $$ \frac{x_{m+1}-x_m}{x_{m}-x_{m-1}}=\frac{\varrho(A_m)}{\varrho(A_{m+1})}\frac{\tilde a_{m}+\tilde d_{m}f_{m-1}}{\tilde a_{m}+\tilde d_{m}\tilde b_{m}\varrho(A_{m})^{-1}}\frac{\tilde d_{m+1}}{\tilde d_m}\frac{\varepsilon_{m}^f}{\varepsilon_{m-1}^f}\rightarrow -1,$$ as $m\rto.$
 So there exists some number $m_2>0$ such that
 $$ \frac{x_{m+1}-x_m}{x_{m}-x_{m-1}}<0\ \text{ for all } m> m_2.$$
Since $\varepsilon_m^f= f_m-\tilde b_{m+1}\varrho(A_{m+1})^{-1}\rightarrow 0$ as $m\rto,$ then by (\ref{xf}), we have $x_{m+1}-x_m\rightarrow0$ as $m\rto.$
We thus come to the conclusion that $x_{m+1}-x_m$ converges to $0$ in an alternating manner as $m\rto.$ Therefore, $$c:=\lim_{m\rto}x_m=x_1+\sum_{m=1}^{\infty} (x_{m+1}-x_m)$$ exists and by (\ref{xul}), we must have $0<c<\infty.$ Theorem \ref{mpr} is proved.\qed
\section{Product of matrices and continued fractions}\label{pc}
Based on Theorem \ref{mpr}, the purpose of this section is to finish the proof of Theorem \ref{s}. Our method is some delicate analysis on continued fractions and their approximants.

\subsection{Product of matrices expressed in terms of the approximants of continued fractions}

Our approach to prove \eqref{sa} is to express $\mb e_1A_k\cdots A_n\mb e_1^t$ in terms of the approximants of some continued fractions. For $1\le k\le n,$ set
 \begin{align}\label{xiy}
   y_{k,n}:=\mathbf e_1 A_k\cdots A_{n}\mathbf e_1^t \text{ and }\xi_{k,n}:=\frac{y_{k+1,n}}{y_{k,n}}.
 \end{align}
Noting that the empty product equals identity, thus $y_{n+1,n}=1.$
 Therefore,
 \begin{align}
\xi_{k,n}^{-1}\cdots \xi_{n,n}^{-1}&=y_{k,n}=\mathbf e_1 A_k\cdots A_{n}\mathbf e_1^t,\label{cpn}\\
\label{axn}\sum_{k=1}^{n+1} \mathbf e_1A_k\cdots A_n\mathbf e_1^t&=\sum_{k=1}^{n+1} \xi_{k,n}^{-1}\cdots \xi_{n,n}^{-1}=\frac{\sum_{k=1}^{n+1} \xi_{1,n}\cdots \xi_{k-1,n}}{\xi_{1,n}\cdots \xi_{n,n}}.
 \end{align}
 \begin{lemma}\label{axc}For $1\le k\le n,$ $\xi_{k,n}$ defined in \eqref{xiy} coincides with the one in \eqref{aprx} with $\beta_k=\tilde b_k^{-1}\tilde d_{k+1}^{-1}$ and $\alpha_k=\tilde a_k\tilde b_k^{-1}\tilde d_{k+1}^{-1}.$
    \end{lemma}
    \proof  Clearly, $\xi_{n,n}=\frac{1}{y_{n,n}}=\frac{1}{\tilde a_n}=\frac{\tilde b_n^{-1}\tilde d_{n+1}^{-1}}{\tilde a_n\tilde b_n^{-1}\tilde d_{n+1}^{-1}}.$ For $1\le k< n,$ note that
\begin{align}\label{ix}
  \xi_{k,n}&=\frac{y_{k+1,n}}{y_{k,n}}=\frac{\mathbf e_1 A_{k+1}\cdots A_{n}\mathbf e_1^t}{\mathbf e_1 A_k\cdots A_{n}\mathbf e_1^t}=\frac{\mathbf e_1 A_{k+1}\cdots A_{n}\mathbf e_1^t}{(\tilde a_k\mathbf e_1+\tilde b_k\mathbf e_2) A_{k+1}\cdots A_{n}\mathbf e_1^t}\\
  &=\frac{1}{\tilde a_k+\tilde b_k\frac{\mathbf e_2 A_{k+1}\cdots A_{n}\mathbf e_1^t}{\mathbf e_1 A_{k+1}\cdots A_{n}\mathbf e_1^t}}=\frac{1}{\tilde a_k+\tilde b_k\tilde d_{k+1}\frac{\mathbf e_1 A_{k+2}\cdots A_{n}\mathbf e_1^t}{\mathbf e_1 A_{k+1}\cdots A_{n}\mathbf e_1^t}}\no\\
  &=\frac{\tilde b_k^{-1}\tilde d_{k+1}^{-1}}{\tilde a_k\tilde b_k^{-1}\tilde d_{k+1}^{-1}+\xi_{k+1,n}}.\no
  \end{align}
  Thus, the lemma can be proved by iterating \eqref{ix}. \qed

In the remainder of this section, we always assume that all conditions of Theorem \ref{s} hold and $\xi_k,\xi_{k,n}, n\ge k\ge1$ are those defined in \eqref{aprx} and \eqref{xic} with $\beta_k=\tilde b_k^{-1}\tilde d_{k+1}^{-1}$ and $\alpha_k=\tilde a_k\tilde b_k^{-1}\tilde d_{k+1}^{-1}.$ Then we have \begin{align*}
  \lim_{k\rto}\beta_k=:\beta=(bd-a\theta)^{-1}\ne 0\text{ and }\lim_{k\rto}\alpha_k=:\alpha=\frac{a+\theta}{bd-a\theta}\ne 0.
\end{align*}
Consequently, it follows from Lemma \ref{ct} that
\begin{align}\label{lxk}
  \lim_{n\rto}\xi_{k,n}=\xi_k, \lim_{k\rto}\xi_k=:\xi=\frac{\alpha}{2}\z(\sqrt{1+4\beta/\alpha^2}-1\y)=\varrho^{-1}>0.
\end{align}
Moreover, by  \eqref{apm1} and \eqref{xiy} we have
\begin{align}\label{po}
 \xi_k>0, \xi_{k,n}>0,\ \forall n\ge k\ge1.
\end{align}
In view of
 \eqref{cpn} and \eqref{axn}, the following proposition  is crucial to prove Theorem \ref{s}.

\begin{proposition}\label{xias}  As $n\rto$, we have
\begin{align}
 &\xi_{1,n}\cdots \xi_{n,n}\sim c\varrho(A_1)^{-1}\cdots\varrho(A_n)^{-1},  \label{xnr}\\
    &\xi_{1,n}\cdots \xi_{n,n}\sim c\xi_{1}\cdots \xi_{n}, \label{xnp}\\
      &\sum_{k=1}^{n+1}\xi_{1,n}\cdots\xi_{k-1,n}\sim c\sum_{k=1}^{n+1}\xi_{1}\cdots \xi_{k-1}. \label{sxnn}
 \end{align}
  \end{proposition}
\proof To begin with, taking \eqref{apm1} and \eqref{cpn} into consideration, applying Theorem \ref{mpr}, we obtain \eqref{xnr}. Since $bd\ne a\theta$ and $\lim_{k\rto}\tilde d_k=b^{-1}(bd-a\theta)\ne 0,$ there exist  $k_0>0$ and $\varepsilon>0$ such that
\begin{align}\label{dk}
 \text{either } \tilde d_k\le -\varepsilon,\ \forall k\ge k_0 \text{ or }  \tilde d_k\ge \varepsilon,\ \forall k\ge k_0.
\end{align}
In view of \eqref{lxk} and \eqref{po}, if we can show that for some $k_0>0,$ $\xi_{k_0,n}\cdots \xi_{n,n}\sim c\xi_{k_0}\cdots \xi_{n}$ and  $\sum_{k=k_0}^{n+1}\xi_{k_0,n}\cdots\xi_{k-1,n}\sim c\sum_{k=k_0}^{n+1}\xi_{k_0}\cdots \xi_{k-1}$ as $n\rto,$ then \eqref{xnp} and \eqref{sxnn} are also true. Therefore,  instead of \eqref{dk}, we assume $
 \exists \varepsilon>0 $  such that either  $\tilde d_k\le -\varepsilon,\ \forall k\ge 1$  or   $\tilde d_k\ge \varepsilon,\ \forall k\ge 1.$ If $ \tilde d_k\ge \varepsilon,\ \forall k\ge 1,$ then Proposition \ref{xias} has been proved in \cite{wy}, see Lemma 10 therein. Thus we need only to deal with the case
 \begin{align}\label{dka}
   \exists \varepsilon>0 \text{ such that }\tilde d_k\le -\varepsilon,\ \forall k\ge 1.
 \end{align}
To prove \eqref{xnp} and \eqref{sxnn}, we need the following two lemmas, whose proof will be postponed to the end of this subsection.

\begin{lemma}\label{mi}Suppose that \eqref{dka} holds and fix $k\ge1.$ Then $\xi_{k,n},n\ge k$ is monotone increasing in $n$ and thus $\xi_{k,n}<\xi_k,n\ge k.$
\end{lemma}

\begin{lemma}\label{dxr} Suppose that \eqref{dka} holds. Then $\exists 0<r<1$ and $ k_1>0$ such that
$
\frac{\xi_{k,n}-\xi_k}{\xi_{k+1,n}-\xi_{k+1}}\le r, \forall n> k>k_1,
$ and consequently
\begin{align*}
 \xi_{k}-\xi_{k,n}\le r^{n-k}(\xi_n-\xi_{n,n}), \forall n\ge k>k_1.
\end{align*}
\end{lemma}
With Lemma \ref{mi} and Lemma \ref{dxr} in hands, we show next that
\begin{align}\label{xnup}
  \exists c>0 \text{ such that } c<\frac{\xi_{1,n}\cdots\xi_{n,n}}{\xi_{1}\cdots\xi_n}<1,n\ge1.
\end{align}
Indeed, applying Lemma \ref{mi}, we get \begin{align}\label{pub}
  \xi_{1,n}\xi_{2,n}\cdots\xi_{n,n}<\xi_1\xi_2\cdots \xi_n,n\ge1.
\end{align}
For a lower bound,
notice that $\xi_{k,n}>0,\forall n\ge k>1$ and $\xi_{n,n}=\tilde a_n^{-1}\rightarrow (a+\theta)^{-1}$ as $n\rto.$ Then by \eqref{lxk}, \eqref{po} and Lemma \ref{mi}, we have \begin{align}\label{sxjx}
  c^{-1}<\xi_n<c \text{ and }c^{-1}<\xi_{n}-\xi_{n,n}<c, \forall n\ge1,
\end{align}
for some number $c>0.$ Therefore,  letting $r$ be the one in Lemma \ref{dxr},   there is a number $k_2>0$ such that
\begin{align*}
  1-\frac{(\xi_n-\xi_{n,n})r^{k_2}}{\inf_{k\ge1}\xi_k}>0, \forall n\ge1.
\end{align*}
With $k_1$ the number in Lemma \ref{dxr}, applying Lemma \ref{dxr}, we have  for $n$ large enough,
\begin{align*}
 \frac{\xi_{k_1+1,n}\cdots\xi_{n,n}}{\xi_{k_1+1}\cdots\xi_n}&=\frac{(\xi_{k_1+1,n}-\xi_{k_1+1})}{\xi_{k_1+1}}\frac{\xi_{k_1+2,n}\cdots\xi_{n,n}}{\xi_{k_1+2}\cdots\xi_n}
 +\frac{\xi_{k_1+2,n}\cdots\xi_{n,n}}{\xi_{k_1+2}\cdots\xi_n}\\
 &\ge \z(1-\frac{(\xi_n-\xi_{n,n})r^{n-k_1-1}}{\xi_{k_1+1}}\y)\frac{\xi_{k_1+2,n}\cdots\xi_{n,n}}{\xi_{k_1+2}\cdots\xi_n},
% \\ &>\frac{\xi_{k_1+2,n}\cdots\xi_{n,n}}{\xi_{k_1+2}\cdots\xi_n}.
\end{align*}
which implies $
 \varliminf_{n\rto}\frac{\xi_{k_1+1,n}\cdots\xi_{n,n}}{\xi_{k_1+1}\cdots\xi_n}\ge  \varliminf_{n\rto} \frac{\xi_{k_1+2,n}\cdots\xi_{n,n}}{\xi_{k_1+2}\cdots\xi_n}.
$
Using the same tricky times and again, we get $$\varliminf_{n\rto}
  \frac{\xi_{k_1+1,n}\cdots\xi_{n,n}}{\xi_{k_1+1}\cdots\xi_n}>\varliminf_{n\rto} \frac{\xi_{n-k_2,n}\cdots\xi_{n,n}}{\xi_{n-k_2}\cdots\xi_n}.
$$  Consequently, since $\lim_{n\rto}\xi_{k,n}=\xi_k$ and $\lim_{n\rto}\frac{\xi_{n-k_2,n}\cdots\xi_{n,n}}{\xi_{n-k_2}\cdots\xi_n}=c$ for some number $c>0,$ we have
\begin{align}\label{plb}
\varliminf_{n\rto}\frac{\xi_{1,n}\cdots\xi_{n,n}}{\xi_{1}\cdots\xi_n}=\varliminf_{n\rto} \frac{\xi_{k_1+1,n}\cdots\xi_{n,n}}{\xi_{k_1+1}\cdots\xi_n}\ge \lim_{n\rto}\frac{\xi_{n-k_2,n}\cdots\xi_{n,n}}{\xi_{n-k_2}\cdots\xi_n}=c>0.
  \end{align}
Taking \eqref{pub} and \eqref{plb} together, we finish the proof of \eqref{xnup}.

Now, we are ready to prove \eqref{xnp}. In view of \eqref{xnr}, it suffices to show
\begin{align}\label{ex}
  \xi_1\cdots \xi_n\sim \varrho(A_1)^{-1}\cdots\varrho(A_n)^{-1}, \text{ as }n\rto.
\end{align}
To this end, write $x_n=\frac{\xi_1\cdots \xi_n}{ \varrho(A_1)^{-1}\cdots\varrho(A_n)^{-1}}, n\ge1.$ Then taking \eqref{xnr} and \eqref{xnup} together, we get
\begin{align}\label{fx}
  C^{-1}<x_n<C,\ \forall n\ge1.
\end{align}
With \eqref{fx} in hand, using Lemma \ref{dxf} and mimicking the proof of the convergence of $x_k$ in the proof of Theorem \ref{mpr}, we can show that $\lim_{n\rto}x_n=c$ for some number $c>0.$ Consequently, \eqref{ex} is proved and so is \eqref{xnp}.

Finally, we turn to prove \eqref{sxnn}. For the case $\varrho<1,$ the proof is the same as the one in \cite[Lemma 10]{wy}. Thus we need only to prove the case $\varrho\ge1.$
Suppose now $\varrho\ge1.$ In order to prove  \eqref{sxnn}, it suffices to show
$
  \lim_{n\rto}\frac{\sum_{k=1}^n\xi_{1,n}\cdots\xi_{k,n}-\sum_{k=1}^n\xi_1\cdots\xi_k}{\sum_{k=1}^n\xi_1\cdots\xi_k}=0
$
which is equivalent to
\begin{align}\label{eqs}
   \lim_{n\rto}\frac{\sum_{j=1}^{n-k_1}\prod_{k=k_1+1}^{k_1+j}\xi_{k,n}}{\sum_{j=1}^{n-k_1}\prod_{k=k_1+1}^{k_1+j}\xi_{k}}=1
\end{align}
since $\lim_{k\rto}\xi_{k,n}=\xi_k>0,$ where $k_1$ is the number in Lemma \ref{dxr}.
To prove \eqref{eqs}, on one hand, note that by Lemma \ref{mi}  we have $\xi_{k,n}<\xi_k,\forall n\ge k\ge1$ which implies
\begin{align}\label{usx}
  \sum_{j=1}^{n-k_1}\prod_{k=k_1+1}^{k_1+j}\xi_{k,n}<\sum_{j=1}^{n-k_1}\prod_{k=k_1+1}^{k_1+j}\xi_{k}, \forall n\ge k_1.
\end{align}
On the other hand, taking the lemmas \ref{mi} and \ref{dxr} into account, we have for $n\ge k_1$
\begin{align*}
  \sum_{j=1}^{n-k_1}&\prod_{k=k_1+1}^{k_1+j}\xi_{k,n}
  =\xi_{k_1+1}\sum_{j=1}^{n-k_1}\prod_{k=k_1+2}^{k_1+j}\xi_{k,n}
  +(\xi_{k_1+1,n}-\xi_{k_1+1})\sum_{j=1}^{n-k_1}\prod_{k=k_1+2}^{k_1+j}\xi_{k,n}  \\
  &\ge\xi_{k_1+1}+\xi_{k_1+1}\sum_{j=2}^{n-k_1}\prod_{k=k_1+2}^{k_1+j}\xi_{k,n}+(\xi_{n,n}-\xi_n)r^{n-(k_1+1)} \sum_{j=2}^{n-k_1}\prod_{k=k_1+1}^{k_1+j}\xi_{k}.
\end{align*}
Iterating the above inequality, we get
\begin{align}\label{lbi}
  \sum_{j=1}^{n-k_1}&\prod_{k=k_1+1}^{k_1+j}\xi_{k,n}>\sum_{j=1}^{n-k_1}\prod_{k=k_1+1}^{k_1+j}\xi_k\\
 &+(\xi_{n,n}-\xi_n)\sum_{i=0}^{n-(k_1+1)}\z(\prod_{k=k_1+1}^{k_1+i}\xi_k\y)r^{n-(k_1+i+1)}\sum_{j=i+1}^{n-k_1}\prod_{k=k_1+i+2}^{k_1+j}\xi_k\no\\
&\quad\quad=:\textrm{(I)}+\textrm{(II)}.\no\end{align}
It follows from \eqref{sxjx} that the absolute value of the second term
(II) equals to
\begin{align} \label{sg}
  (\xi_n&-\xi_{n,n})\sum_{i=0}^{n-(k_1+1)}\frac{r^{n-(k_1+i+1)}}{\xi_{k_1+i+1}}\z(\prod_{k=k_1+1}^{k_1+i+1}\xi_k\y)\sum_{j=i+1}^{n-k_1}\prod_{k=k_1+i+2}^{k_1+j}\xi_k\\
&\le c \sum_{i=0}^{n-(k_1+1)}{r^{n-(k_1+i+1)}}\z(\prod_{k=k_1+1}^{k_1+i+1}\xi_k\y)\sum_{j=i+1}^{n-k_1}\prod_{k=k_1+i+2}^{k_1+j}\xi_k.\no
\end{align}
Clearly, for any $0\le i\le n-(k_1+1),$ \begin{align}\label{fla}
 \frac{\z(\prod_{k=k_1+1}^{k_1+i+1}\xi_k\y)\sum_{j=i+1}^{n-k_1}\prod_{k=k_1+i+2}^{k_1+j}\xi_k}{ \sum_{j=1}^{n-k_1}\prod_{k=k_1+1}^{k_1+j}\xi_k}\le1\end{align}
and an application of \cite[Lemma 8]{wy} yields that
\begin{align}\label{flb}
  \lim_{n\rto}\frac{\prod_{k=k_1+1}^{n-i}\xi_k}{\sum_{j=1}^{n-k_1}\prod_{k=k_1+1}^{k_1+j}\xi_k}=0, \forall i\ge1,
\end{align} since $\xi_k>0,\forall k\ge1$ and $\lim_{k\rto}\xi_k=\varrho^{-1}\le 1.$

Let $\varepsilon>0$ be an arbitrary number and fix $k_3$ such that $r^{k_3}<\varepsilon.$
Then using \eqref{fla} and \eqref{flb}, we have from \eqref{sg} that
\begin{align*}
 \varlimsup_{n\rto}& \frac{|\textrm{(II)}|}{\sum_{j=1}^{n-k_1}\prod_{k=k_1+1}^{k_1+j}\xi_k}\le c\varlimsup_{n\rto}\sum_{i=0}^{n-(k_1+1)-k_3}r^{n-(k_1+i+1)}\\
 &+c\varlimsup_{n\rto}\sum_{i=n-k_1-k_3}^{n-(k_1+1)}r^{n-(k_1+i+1)}
 \frac{\z(\prod_{k=k_1+1}^{k_1+i+1}\xi_k\y)\sum_{j=i+1}^{n-k_1}\prod_{k=k_1+i+2}^{k_1+j}\xi_k}{\sum_{j=1}^{n-k_1}\prod_{k=k_1+1}^{k_1+j}\xi_k}\no\\
 &=c\varlimsup_{n\rto}\sum_{i=0}^{n-(k_1+1)-k_3}r^{n-(k_1+i+1)} =cr^{k_3}/(1-r)\no\\
 &\le c\varepsilon/(1-r).\no
\end{align*}
Since $\varepsilon>0$ is arbitrary, $  \varlimsup_{n\rto} \frac{|\textrm{(II)}|}{\sum_{j=1}^{n-k_1}\prod_{k=k_1+1}^{k_1+j}\xi_k}=0,$ which together  with \eqref{lbi} implies that
\begin{align}\label{lsx}
  \varliminf_{n\rto}\frac{\sum_{j=1}^{n-k_1}\prod_{k=k_1+1}^{k_1+j}\xi_{k,n}}{\sum_{j=1}^{n-k_1}\prod_{k=k_1+1}^{k_1+j}\xi_k}=1.
\end{align}
Putting \eqref{usx} and \eqref{lsx} together, we can infer that \eqref{eqs} is true and thus Proposition \ref{xias} is proved. \qed

To end this subsection, we give the proof of Lemma \ref{mi} and Lemma \ref{dxr}.

\noindent {\it Proof of Lemma \ref{mi}.}
 We see from \eqref{apm1} and \eqref{xiy} that $\xi_{k,n}>0$ for all $n\ge k\ge1.$ Since $\forall k\ge1,$ $\tilde a_k>0, \tilde b_k>0$ and $\tilde d_k<-\varepsilon<0$ for some $\varepsilon>0,$ we always have $\alpha_k<0,\beta_k<0,\forall k\ge1.$
 Now fix $k\ge1$ and for $n\ge k,$ write
\begin{align}\label{cd}
C_{k,n}&=\alpha_nC_{k,n-1}+\beta_{n}C_{k,n-2},\
D_{k,n}=\alpha_nD_{k,n-1}+\beta_{n}D_{k,n-2}
\end{align}
with initial values \begin{align}
   C_{k-2}=D_{k-1}=1, C_{k-1}=D_{k-2}=0.\label{cdi}
 \end{align}
 Then by Euler-Minding formula (see \cite{lw}, (1.2.14) on page 7) we obtain
\begin{align*}
  \xi_{k,n}=\frac{C_{k,n}}{D_{k,n}}=-\sum_{j=k}^{n}\frac{\prod_{i=k}^j(-\beta_i)}{D_{k,j}D_{k,j-1}}
\end{align*}
which leads to
\begin{align}\label{dxikn}
  \xi_{k,n+1}-\xi_{k,n}=-\frac{\prod_{i=k}^{n+1}(-\beta_j)}{D_{k,n+1}D_{k,n}},n\ge k.
\end{align}
By induction, we have from \eqref{cd} and \eqref{cdi} that
\begin{align*}
&\left(
    \begin{array}{cc}
      D_{k,n} & D_{k,n-1} \\
      C_{k,n} & C_{k,n-1} \\
    \end{array}
  \right)=\left(
            \begin{array}{cc}
              \alpha_k & 1 \\
              \beta_k & 0 \\
            \end{array}
          \right)\cdots\left(
            \begin{array}{cc}
              \alpha_n & 1 \\
              \beta_n & 0 \\
            \end{array}
          \right)\\
  &\quad\quad={\z(\prod_{i=k}^n\tilde b_i\tilde d_{i+1}\y)^{-1}}{\left(
            \begin{array}{cc}
              \tilde a_k & \tilde b_k\tilde d_{k+1} \\
              1 & 0 \\
            \end{array}
          \right)\cdots\left(
            \begin{array}{cc}
              \tilde a_n & \tilde b_{n}\tilde d_{n+1}  \\
              1 & 0 \\
            \end{array}
          \right)}\no\\
 &\quad\quad={\z(\prod_{i=k}^n\tilde b_i\tilde d_{i+1}\y)^{-1}}{\left(
            \begin{array}{cc}
              1 & 0 \\
              0 & \tilde d_k^{-1} \\
            \end{array}
          \right)A_k\cdots A_n\left(
            \begin{array}{cc}
              1 & 0 \\
              0 & \tilde d_{n+1} \\
            \end{array}
          \right)}\no
\end{align*} where for the third equality we use the facts
\begin{align*}
  &\left(
            \begin{array}{cc}
              \tilde a_i & \tilde b_i\tilde d_{i+1} \\
              1 & 0 \\
            \end{array}
          \right)=\left(
            \begin{array}{cc}
              \tilde a_i & \tilde b_i \\
              1 & 0 \\
            \end{array}
          \right)\left(
            \begin{array}{cc}
              1 & 0  \\
              0 & \tilde d_{i+1} \\
            \end{array}
          \right),\\
         & A_i=\left(
            \begin{array}{cc}
              \tilde a_i & \tilde b_i \\
              \tilde d_{i} & 0 \\
            \end{array}
          \right)=\left(
            \begin{array}{cc}
              1 & 0  \\
              0 &\tilde d_{i} \\
            \end{array}
          \right) \left(
            \begin{array}{cc}
              \tilde a_i & \tilde b_i \\
              1 & 0 \\
            \end{array}
          \right),i\ge1.
\end{align*}
As a consequence, taking \eqref{apm1} and \eqref{apm2} into account we get
\begin{align}\label{dd}
  \frac{D_{k,n-1}}{D_{k,n}}=\tilde d_{n+1}\frac{\mb e_1A_k\cdots A_n\mb e_2}{\mb e_1A_k\cdots A_n\mb e_1}<0,\forall n\ge k.
\end{align}
Putting \eqref{dxikn} and \eqref{dd} together, we conclude that
\begin{align}\no
  \xi_{k,n+1}-\xi_{k,n}>0, \forall n\ge k\ge1,
\end{align}
that is to say, $\xi_{k,n},n\ge k$ is monotone increasing in $n.$ Lemma \ref{mi} is proved. \qed

%%%%%%%%%%%%%%%%%%%%%%%%%%%%%%%%%%%%%%%%%%%%%%%%%%%%%%%%%%%%%%%%%%%%%%%%%%%%%%%%%

\noindent{\it Proof of Lemma \ref{dxr}}.
 Taking Lemma \ref{mi} into account, since $\xi_{k,n}>0,$ $\forall n\ge k\ge1,$ then by some easy computation, we have
 from \eqref{aprx} and \eqref{xic} that
 \begin{align*}
 0<\frac{\xi_{k,n}-\xi_k}{\xi_{k+1,n}-\xi_{k+1}}=-\frac{\xi_{k,n}}{\alpha_k+\xi_{k+1}}<-\frac{\xi_{k}}{\alpha_k+\xi_{k+1}}.
 \end{align*}
But
 $-\frac{\xi_{k}}{\alpha_k+\xi_{k+1}}\rightarrow-\frac{\xi}{\alpha+\xi}=\frac{\alpha+\sqrt{\alpha^2+4\beta}}{\alpha-\sqrt{\alpha^2+4\beta}}<1$ as $k\rto.$
 As a result,  for some proper number $0<r<1,$ $\exists k_1>0$ such that
 $-\frac{\xi_{k}}{\alpha_k+\xi_{k+1}}<r, \forall k> k_1,$ which finishes the proof of the lemma. \qed

 \subsection{Proof of Theorem \ref{s}}
 Theorem \ref{s} is a direct consequence of Proposition \ref{xias}. In fact, it follows from \eqref{axn} and Proposition \ref{xias} that
 \begin{align*}
   \sum_{k=1}^{n+1} \mathbf e_1A_k\cdots A_n\mathbf e_1^t=\frac{\sum_{k=1}^{n+1} \xi_{1,n}\cdots \xi_{k-1,n}}{\xi_{1,n}\cdots \xi_{n,n}}\sim c\frac{\sum_{k=1}^{n+1}\xi_{1}\cdots \xi_{k-1}}{\varrho(A_1)^{-1}\cdots \varrho(A_n)^{-1}}, \text{ as }n\rto.
 \end{align*}
 But by \eqref{xnr} and \eqref{xnp}, we obtain
 \begin{align*}
   \xi_{1}\cdots \xi_{n} \sim c\varrho(A_1)^{-1}\cdots \varrho(A_n)^{-1}, \text{ as }n\rto
 \end{align*}
 which implies that
 \begin{align*}
   \sum_{k=1}^{n+1}\xi_{1}\cdots \xi_{k-1}\sim c\sum_{k=1}^{n+1}\varrho(A_1)^{-1}\cdots \varrho(A_{k-1})^{-1}, \text{ as }n\rto.
 \end{align*}
 As a result, we have
 \begin{align*}
   \sum_{k=1}^{n+1} \mathbf e_1A_k\cdots A_n\mathbf e_1^t\sim c\frac{\sum_{k=1}^{n+1}\varrho(A_1)^{-1}\cdots \varrho(A_{k-1})^{-1}}{\varrho(A_1)^{-1}\cdots \varrho(A_n)^{-1}}=c\sum_{k=1}^{n+1} \varrho(A_k)\cdots \varrho(A_{n})
 \end{align*}
 as $n\rto.$ Thus Theorem \ref{s} is proved. \qed
%%%%%%%%%%%%%%%%%%%%%%%%%%%%%%%%%%%%%%%%%%%%%%%%%%%%%%%%%%%%%%%%%%%%%%%%%%%%%%%%%%%%%%%%%%%%%%%%%%5

\section{Proof of Theorem \ref{pnec}}\label{pm}

For $n\ge 0$ and $\mathbf s=(s_1,s_2)^t\in [0,1]^2,$  let $$F_n^{(i)}(\mathbf s)\equiv E(\mathbf s^{Z_n}|Z_0=\mathbf e_i):=E( s_1^{Z_{n,1}}s_2^{Z_{n,2}}|Z_0=\mathbf e_i),i=1,2$$
  and set $\mathbf F_n(\mathbf s)=(F_n^{(1)}(\mathbf s),F_n^{(2)}(\mathbf s))^t.$
It follows by induction (see Dyakonova \cite[Lemma 1]{dy99}) that for $n\ge0,$
\begin{align*}
  \mathbf F_{n}(\mathbf{s})=\mathbf{f}_{1}(\mathbf{f}_{2}(\cdots\mathbf{f}_{n}(\mathbf{s})\cdots))
  =\mathbf1-\frac{\prod_{k=1}^{n}M_{k}(\mathbf 1-\mathbf{s})}{1+\sum_{k=1}^{n}\gamma_k\prod_{i=k+1}^{n}M_{i}(\mathbf1-\mathbf{s})}
\end{align*}
which leads to
\begin{align*}
  \mathbf F_{n}(\mathbf{0})=\mathbf 1- \frac{\prod_{k=1}^{n}M_{k}\mathbf 1}{1+\sum_{k=1}^{n}\gamma_k\prod_{i=k+1}^{n}M_{i}\mathbf1}.
\end{align*}

For the asymptotics of $P(\nu>n|Z_0=\mb e_1)$ and $P(\nu=n|Z_0=\mb e_1),$ we refer the reader to \cite{wy}. Here we treat only $P(\nu>n|Z_0=\mb e_2)$ and $P(\nu=n|Z_0=\mb e_2).$  Taking \eqref{mg} intro account, for $n\ge1$ we have
\begin{align}\label{pngn}
  P(\nu>n|Z_0=\mb e_2)&=1-F_n^{(2)}(\mb 0)=\frac{\mb e_2\prod_{k=1}^{n}M_{k}\mathbf 1}{\sum_{k=1}^{n+1}\mb e_1\prod_{i=k}^{n}M_{i}\mathbf1},
\end{align}
and consequently \begin{align}\label{pnen}
  P(&\nu=n|Z_0=\mb e_2)=\frac{\mb e_2\prod_{k=1}^{n-1}M_{k}\mathbf 1}{\sum_{k=1}^{n}\mb e_1\prod_{i=k}^{n-1}M_{i}\mathbf1}-\frac{\mb e_2\prod_{k=1}^{n}M_{k}\mathbf 1}{\sum_{k=1}^{n+1}\mb e_1\prod_{i=k}^{n}M_{i}\mathbf1}.
\end{align}

With matrices  $A_i,i\ge1$ the ones defined in \eqref{dta},
   using  \eqref{am}, by some very careful computation, we have from \eqref{pngn} and \eqref{pnen} that
 \begin{align}\label{etan}
     P(\nu>n|Z_0=\mb e_2)&=\frac{(\theta_1/b_1,1)\prod_{k=1}^{n}A_k(1,\lambda_{n+1})^t}{\sum_{k=1}^{n+1}\mb e_1\prod_{i=k}^{n}A_i(1,\lambda_{n+1})^t},\\
\label{pn}
  P(\nu=n|Z_0=\mb e_2)&=\frac{1}{\sum_{k=1}^{n+1}\mb e_1\prod_{i=k}^{n}A_i(1,\lambda_{n+1})^t}\frac{\mb e_1\prod_{k=1}^{n-1}A_{k}\mathbf e_1^t}{\sum_{k=1}^{n}\mb e_1\prod_{i=k}^{n-1}A_i(1,\lambda_{n})^t}\\
  &\quad\quad\times\z(\frac{\theta_1}{b_1}G_{n-1,1}+\frac{\mb e_2\prod_{k=1}^{n-1}A_{k}\mathbf e_1^t}{\mb e_1\prod_{k=1}^{n-1}A_{k}\mathbf e_1^t}G_{n-1,2}\y)\no
  \end{align}
 where  for $n\ge1,$  $\lambda_n\equiv 1-\frac{\theta_n}{b_n},$
\begin{align}\label{dg}
  G_{n-1,1}\equiv \frac{\mb e_1\prod_{k=1}^{n-1}A_k(1,\lambda_{n})^t\sum_{k=1}^{n+1}\mb e_1\prod_{i=k}^{n}A_i(1,\lambda_{n+1})^t}{\mb e_1\prod_{k=1}^{n-1}A_{k}\mathbf e_1^t}\\
  -\frac{\mb e_1\prod_{k=1}^{n}A_{k}(1,\lambda_{n+1})^t\sum_{k=1}^{n}\mb e_1\prod_{i=k}^{n-1}A_i(1,\lambda_{n})^t}{\mb e_1\prod_{k=1}^{n-1}A_{k}\mathbf e_1^t}\no
  \end{align}
and
\begin{align}\label{dg2}
  G_{n-1,2}\equiv \frac{\mb e_2\prod_{k=1}^{n-1}A_k(1,\lambda_{n})^t\sum_{k=1}^{n+1}\mb e_1\prod_{i=k}^{n}A_i(1,\lambda_{n+1})^t}{\mb e_1\prod_{k=1}^{n-1}A_{k}\mathbf e_1^t}\\
  -\frac{\mb e_2\prod_{k=1}^{n}A_{k}(1,\lambda_{n+1})^t\sum_{k=1}^{n}\mb e_1\prod_{i=k}^{n-1}A_i(1,\lambda_{n})^t}{\mb e_1\prod_{k=1}^{n-1}A_{k}\mathbf e_1^t}\no.
  \end{align}
 $G_{n,1}$ and $G_{n,2}$ defined in \eqref{dg} and \eqref{dg2} look very complicated. Next lemma shows that they converge to the same limit.
 \begin{lemma}\label{gl}
  Suppose that condition (B1) holds and $|\varrho_1|<1.$  Then, there exists a number $G\ge 0$ such that \begin{align}\label{lg}
  \lim_{n\rto}G_{n,2}=G=\lim_{n\rto}G_{n,1}.\end{align}  Moreover,  $G\ne0$ if and only if  $\varrho_1\ne\frac{1}{2}\z(a+b+1-\sqrt{(a+b+1)^2+4\frac{bd- a\theta}{\theta-b}}\y).$
\end{lemma}
\proof The proof of the second equality in \eqref{lg} can be found in \cite{wy}, so we check here only the first one. To this end, set \begin{align*}
  f_n\equiv \frac{\mb e_2\prod_{k=1}^{n}A_{k}\mathbf e_2^t}{\mb e_2\prod_{k=1}^{n}A_{k}\mathbf e_1^t} \text{ and }H_n\equiv\sum_{k=1}^{n}\mb e_1\prod_{i=k}^{n}A_{i}(f_n\mathbf e_1^t-\mb e_2^t),n\ge1.
\end{align*}
Clearly, we have
$
  f_1=0, H_1=-\tilde b_1=-b_1
$
and  by some subtle computation,
\begin{align}\label{ghf}
  G_{n-1,2}=1&+(\tilde b_n\lambda_n\lambda_{n+1} +\tilde a_n\lambda_{n}-\tilde d_n)H_{n-1}\\
  &+(\tilde b_n\lambda_n\lambda_{n+1} +\tilde a_n\lambda_{n}-\tilde d_n+\lambda_{n})f_{n-1}.\no
\end{align}
Note that
\begin{align}
  \label{fr}f_n=\frac{\mb e_2A_1\cdots A_{n}\mathbf e_2^t}{\mb e_2A_{1}\cdots A_n\mathbf e_1^t}=\frac{\tilde b_n\mb e_2A_1\cdots A_{n-1}\mb e_1^t}{\mb e_2A_{1}\cdots A_{n-1}(\tilde a_n\mathbf e_1^t+\tilde d_n\mb e_2^t)}=\frac{\tilde b_n}{\tilde a_n+\tilde d_nf_{n-1}},
\end{align}
which leads to
\begin{align}
  \tilde d_nf_nf_{n-1}=\tilde b_n-\tilde a_nf_{n},n\ge2.\no
\end{align}
Consequently, for $n\ge2,$
\begin{align}
  H_n&=\sum_{k=1}^{n}\mb e_1\prod_{i=k}^{n}A_{i}(f_n\mathbf e_1^t-\mb e_2^t)\label{hr}\\
  &=\tilde a_nf_n-\tilde b_n+\sum_{k=1}^{n-1}\mb e_1\prod_{i=k}^{n-1}A_{i}((\tilde a_nf_n-\tilde b_n)\mathbf e_1^t+\tilde d_nf_n\mb e_2^t)\no\\
  &=- \tilde d_nf_nf_{n-1} -\tilde d_nf_n\sum_{k=1}^{n-1}\mb e_1\prod_{i=k}^{n-1}A_{i}(f_{n-1}\mathbf e_1^t-\mb e_2^t)\no\\
  &=- \tilde d_nf_nf_{n-1} -\tilde d_nf_n H_{n-1}.\no
\end{align}
Since $f_1=0$ and $H_1=-b_1,$ iterating \eqref{hr}, we get
\begin{align}\label{h}
  H_n=(-1)^n b_1\tilde d_2f_2\cdots \tilde d_nf_n+\sum_{k=1}^{n-1}(-1)^{n-k}f_k\tilde d_{k+1}f_{k+1}\cdots \tilde d_nf_n,n\ge2
\end{align}
and iterating \eqref{fr}, we get
\begin{align*}
  f_n=\frac{\tilde b_n\tilde d_n^{-1}}{\tilde a_n\tilde d_n^{-1}}\begin{array}{c}
                                \\
                               +
                             \end{array}\frac{\tilde b_{n-1}\tilde d_{n-1}^{-1}}{ \tilde a_{n-1}\tilde d_{n-1}^{-1}}\begin{array}{c}
                                \\
                               +\cdots+
                             \end{array}\frac{\tilde b_2\tilde d_2^{-1}}{\tilde a_2\tilde d_2^{-1}},n\ge2.
\end{align*}
Noticing that $\lim_{n\rto}\tilde b_n\tilde d_n^{-1}=b^2(bd-a\theta)^{-1}\ne0$ and $\lim_{n\rto}\tilde a_n\tilde d_n^{-1}=b(a+\theta)(bd-a\theta)^{-1},$ thus applying Lemma \ref{ct}, we get
\begin{align}\label{lf}
 \lim_{n\rto}f_n=-\frac{b\varrho_1}{bd-a\theta}
\end{align}
and thus
$
  \lim_{n\rto}\tilde d_nf_n=-\varrho_1.
  $
  Since   $|\varrho_1|<1,$ we conclude from  \eqref{h} that
  \begin{align}\label{lh}
    \lim_{n\rto} H_n=-\frac{b\varrho_1^2}{(bd-a\theta)(1-\varrho_1)}.
  \end{align}
Letting $n\rto$ in \eqref{ghf}, owing to \eqref{lf} and \eqref{lh}, we get
\begin{align*}
  \lim_{n\rto}G_{n-1,2}=G:=\frac{(b-\theta)\varrho_1^2-(b-\theta)(a+b+1)\varrho_1+bd-a\theta }{(bd-a\theta)(1-\varrho_1)}.
\end{align*}
Consequently, the first equality in \eqref{lg} is true.
Moreover, since   $|\varrho_1|<1,$ then by the definition of $G,$ it is easy to see that
$$G\ne 0\text{ if and only if } \varrho_1\ne\frac{1}{2}\z(a+b+1-\sqrt{(a+b+1)^2+4\frac{bd- a\theta}{\theta-b}}\y).$$

What's left for us to do is  to check the  nonnegativity of the number $G.$
For this purpose,  applying  \eqref{lg} and taking \eqref{apm1}-\eqref{apm21} and Theorem \ref{mpr} into account, we have
\begin{align}\label{cg}
  \lim_{n\rto}&\z(\frac{\theta_1}{b_1}G_{n-1,1}+\frac{\mb e_2\prod_{k=1}^{n-1}A_{k}\mathbf e_1^t}{\mb e_1\prod_{k=1}^{n-1}A_{k}\mathbf e_1^t}G_{n-1,2}\y)\\&=G\lim_{n\rto}\z(\frac{\theta_1}{b_1}+\frac{\mb e_2\prod_{k=1}^{n-1}A_{k}\mathbf e_1^t}{\mb e_1\prod_{k=1}^{n-1}A_{k}\mathbf e_1^t}\y)=GL\no
\end{align} exists, where using again \eqref{apm1}-\eqref{apm21} and Theorem \ref{mpr}, we see that
\begin{align*}
  L:&=\lim_{n\rto}\z(\frac{\theta_1}{b_1}+\frac{\mb e_2\prod_{k=1}^{n-1}A_{k}\mathbf e_1^t}{\mb e_1\prod_{k=1}^{n-1}A_{k}\mathbf e_1^t}\y)=\lim_{n\rto}\frac{\mb e_2\prod_{k=1}^{n-1}M_k(1,\theta_n/b_n)^t}{\mb e_1\prod_{k=1}^{n-1}M_k(1,\theta_n/b_n)^t}\no\\
  &=\lim_{n\rto}\frac{\mb e_2\prod_{k=1}^{n-1}M_k(1,\theta_n/b_n)^t}{\mb e_1\prod_{k=1}^{n-1}A_{k}\mathbf e_1^t}\ge \lim_{n\rto}\frac{ d_1d_{n-1}\mb e_1\prod_{k=2}^{n-2}M_k\mb e_2^t}{\mb e_1\prod_{k=1}^{n-1}A_{k}\mathbf e_1^t}\no\\
   &=\lim_{n\rto}\frac{ d_1d_{n-1}\mb e_1\prod_{k=2}^{n-2}A_k\mb e_2^t}{\mb e_1\prod_{k=1}^{n-1}A_{k}\mathbf e_1^t}>0.\no
\end{align*}
 But in view of \eqref{pn}, since  $\sum_{k=1}^{n}\mb e_1\prod_{i=k}^{n-1}A_i(1,\lambda_{n})^t=\sum_{k=1}^{n}\mb e_1\prod_{i=k}^{n-1}M_{i}\mathbf1>0$ and $\mb e_1\prod_{k=1}^{n-1}A_{k}\mathbf e_1^t=\mb e_1\prod_{k=1}^{n-1}M_{k}\mathbf (1,\theta_{k+1}/b_{k+1})^t>0,$ we must have $GL\ge 0$ and consequently, $G\ge 0.$
The lemma is proved. \qed

Next we continue with the proof of Theorem \ref{pnec}. For $n\ge1,$ write   $$S_n:=\frac{\sum_{k=1}^{n+1}\varrho(A_1)^{-1}\cdots\varrho(A_{k-1})^{-1}}{\varrho(A_1)^{-1}\cdots\varrho(A_n)^{-1}}.
 $$
 We need in addition the following  lemma.
 \begin{lemma} \label{smr} Under the conditions of Theorem \ref{pnec}, we have
 \begin{align}\label{epar}
   &\mb e_1\prod_{i=1}^{n}A_n\mb e_1^t\sim c\varrho(A_1)\cdots\varrho(A_{n}),\\
   &\label{epab}\mb e_1\prod_{i=1}^{n}A_n\mb e_2^t\sim c\varrho(A_1)\cdots\varrho(A_{n}),\\
   &\label{epac}(\theta_1/b_1,1)\prod_{k=1}^{n}A_k(1,\lambda_{n+1})^t\sim c \varrho(A_1)\cdots\varrho(A_{n}),
 \end{align} as $n\rto$ and for some number $0<\phi_1,\phi_2<\infty,$
   \begin{align}\label{asal}
    \lim_{n\rto}&\frac{\sum_{k=1}^{n+1}\mb e_1\prod_{i=k}^{n}A_i(1,\lambda_{n+1})^t}{S_n}=\phi_1 \text{ and } \lim_{n\rto}\frac{S_{n-1}}{S_n}=\phi_2.
  \end{align}
 \end{lemma}
 \proof Applying Theorem \ref{mpr}, since  $\mb e_1\prod_{i=1}^{n}A_n\mb e_i^t>0,i\in\{1,2\}$ by \eqref{apm1} and \eqref{apm2}, we get \eqref{epar} and \eqref{epab}. Also, we infer from Theorem \ref{mpr} that
  $$\lim_{n\rto}\frac{(\theta_1/b_1,1)\prod_{k=1}^{n}A_k(1,\lambda_{n+1})^t}{\varrho(A_1)\cdots\varrho(A_{n})}$$ exists. But using \eqref{epab}, we have
  \begin{align}
    \lim_{n\rto}&\frac{ (\theta_1/b_1,1)\prod_{k=1}^{n}A_k(1,\lambda_{n+1})^t}{\varrho(A_1)\cdots\varrho(A_{n})}
    =\lim_{n\rto}\frac{ \mb e_2\prod_{k=1}^{n}M_k(1,1)^t}{\varrho(A_1)\cdots\varrho(A_{n})}\no\\
    & \ge \lim_{n\rto}\frac{ d_1\mb e_1\prod_{k=2}^{n}M_k\mb e_2^t}{\varrho(A_1)\cdots\varrho(A_{n})}=    \lim_{n\rto}\frac{ d_1\mb e_1\prod_{k=2}^{n}A_k\mb e_2^t}{\varrho(A_1)\cdots\varrho(A_{n})}>0.\no
  \end{align}
So  \eqref{epac} is proved true. Finally, with Theorem \ref{s} in hand, the proof of \eqref{asal} is similar to the one of \cite[Lemma 4]{wy} and we will not repeat it here.\qed

 Now we are ready to finish the proof of Theorem \ref{pnec}.  Note that from \eqref{pn} and the lemmas \ref{bma} and \ref{smr}, we have
  \begin{align}
     \no P(&\nu=n|Z_0=\mb e_2)\sim c\frac{\varrho(A_1)\cdots\varrho(A_n)}{S_n^2}\frac{S_n}{S_{n-1}}\z(\frac{\theta_1}{b_1}G_{n-1,1}+\frac{\mb e_2\prod_{k=1}^{n-1}A_{k}\mathbf e_1^t}{\mb e_1\prod_{k=1}^{n-1}A_{k}\mathbf e_1^t}G_{n-1,2}\y)\no\\
     &\sim c\frac{\varrho(M_1)^{-1}\cdots \varrho(M_n)^{-1}}{\z(\sum_{k=1}^{n+1} \varrho(M_1)^{-1}\cdots \varrho(M_{k-1})^{-1}\y)^2}\z(\frac{\theta_1}{b_1}G_{n-1,1}+\frac{\mb e_2\prod_{k=1}^{n-1}A_{k}\mathbf e_1^t}{\mb e_1\prod_{k=1}^{n-1}A_{k}\mathbf e_1^t}G_{n-1,2}\y).\no
  \end{align}
    Then it follows from \eqref{cg} that
    \begin{align}\no
\lim_{n\rto}\frac{P(\nu=n|Z_0=\mb e_2)\z(\sum_{k=1}^{n+1} \varrho(M_1)^{-1}\cdots \varrho(M_{k-1})^{-1}\y)^2}{\varrho(M_1)^{-1}\cdots \varrho(M_n)^{-1}}=GL\ge0,
    \end{align}
    where $L>0$ is a proper constant.  But by Lemma \ref{gl}, $G>0$ if and only if $\varrho_1\ne\frac{1}{2}\z(a+b+1-\sqrt{(a+b+1)^2+4\frac{bd- a\theta}{\theta-b}}\y).$
    We thus conclude that \eqref{pnse} and \eqref{pnso} hold.

    Finally, we turn to prove \eqref{pns}. Note that by \eqref{etan}, \eqref{epac} and \eqref{asal}, we obtain
    \begin{align}\no
      P(\nu>n|Z_0=\mb e_2)\sim c\frac{\varrho(A_1) \cdots\varrho(A_n)}{S_n}=\frac{c}{\sum_{k=1}^{n+1}\varrho(A_1)^{-1}\cdots\varrho(A_{k-1})^{-1}},
    \end{align}
   as $n\rto.$ As a result, applying Lemma \ref{bma}, we get \eqref{pns}. Theorem \ref{pnec} is proved.
    \qed

\vspace{.5cm}

\noindent{{\bf \Large Acknowledgements:}} %The authors would like to thank Prof.  Hong, W.M. for introducing to us the basics of BPVE  and  Prof. Vatutin, V. for some discussions on the distribution of the extinction time when writing the paper. Finally the authors were in debt to two
%referees who read the paper carefully and gave very good suggestions which help to improve the paper to a large extent.
This project is partially supported by National
Natural Science Foundation of China (Grant No. 12071003).

%\newpage

\end{document}